\documentstyle[12pt,newlfont,amssymb]{amsart}

\newcommand{\nc}{\newcommand}

\def\binome#1#2#3{\setlength{\arraycolsep}{0.2truemm}{\small\left[
\begin{array}{c} {#2}\\ {#3}\end{array}\right]}_{#1}}

\nc{\one}{\mbox{\bf 1}}
\nc{\invtensor}{\underset{\leftarrow}{\otimes}}
\nc{\ad}{\operatorname{ad}}
\nc{\wad}{\widetilde{\operatorname{ad}}}
\nc{\adp}{\operatorname{ad}_{\bullet}}
\nc{\badp}{\overline{\operatorname{ad}}_{\bullet}}
\nc{\bad}{\overline{\operatorname{ad}}}
\nc{\rk}{\operatorname{rank}}
\nc{\corank}{\operatorname{corank}}
\nc{\Sym}{\operatorname{Sym}}
\nc{\sym}{\operatorname{sym}}
\nc{\id}{\operatorname{id}}
\nc{\htt}{\operatorname{ht}}
\nc{\Ker}{\operatorname{Ker}}
\nc{\im}{\operatorname{Im}}
\nc{\osp}{\operatorname{osp}}
\nc{\oo}{\operatorname{o}}
\nc{\ssp}{\operatorname{sp}}
\nc{\sgn}{\operatorname{sgn}}
\nc{\F}{\operatorname{F}}
\nc{\bFp}{\overline{\operatorname{F}}_{\bullet}}
\nc{\intl}{\operatorname{int}}
\nc{\out}{\operatorname{out}}
\nc{\Tor}{\operatorname{Tor}}
\nc{\Hom}{\operatorname{Hom}}
\nc{\End}{\operatorname{End}}
\nc{\supp}{\operatorname{supp}}

\nc{\Ann}{\operatorname{Ann}}
\nc{\Ind}{\operatorname{Ind}}
\nc{\Card}{\operatorname{Card}}

\nc{\wt}{\operatorname{wt}}
\nc{\ch}{\operatorname{ch}}
\nc{\Stab}{\operatorname{Stab}}
\nc{\Sch}{{\cal S}\mbox{\em ch}}

\nc{\Spec}{\operatorname{Spec}}
\nc{\Prim}{\operatorname{Prim}}
\nc{\Aut}{\operatorname{Aut}}
\nc{\Fract}{\operatorname{Fract}}
\nc{\gr}{\operatorname{gr}}
\nc{\wdM}{\widetilde{M}}
\nc{\wdV}{\widetilde{V}}

\nc{\xa}{y^{-1}_{-\alpha}}
\nc{\xb}{y^{-1}_{-\alpha}y_{-\alpha-\beta}}
\nc{\xc}{\ol{y}_{-\alpha-\beta}}
\nc{\xd}{y_{-\gamma}}
\nc{\xe}{y^{-1}_{-\alpha}y_{-\alpha-\beta-\gamma}}
\nc{\xf}{\ol{y}_{-\alpha-\beta-\gamma}}
\nc{\xg}{y_{-\beta}}
\nc{\xh}{y_{-\beta-\gamma}}
\nc{\xy}{y_{-\alpha-2\beta-\gamma}}
\nc{\xj}{\ol{y}_{-\beta-\gamma}}

\nc{\Dglie}{\operatorname{{\cal D}glie}}
\nc{\Dgcoalg}{\operatorname{{\cal D}gcoalg}}
\nc{\Homcoalg}{\operatorname{{\cal H}omcoalg}}
\nc{\Homlie}{\operatorname{{\cal H}omlie}}
\nc{\Dgmod}{\operatorname{{\cal D}gmod}}
\nc{\Dgcomod}{\operatorname{{\cal D}gcomod}}

\nc{\pa}{\partial}
\nc{\co}{\cal O}

\nc{\Sg}{{\cal S}({\frak g})}
\nc{\wU}{\widehat{\cal U}}
\nc{\bU}{\overline{\cal U}}
\nc{\wF}{\widetilde{\operatorname{F}}}
\nc{\bF}{\overline{\operatorname{F}}}
\nc{\Ug}{{\cal U}({\frak g})}
\nc{\Uq}{{\cal U}_q}
\nc{\Uhg}{{\cal U}_h({\frak g})}
\nc{\Uqg}{{\cal U}_q({\frak g})}
\nc{\Urs}{{\cal U}_{h,r_s}({\frak g})}
\nc{\Ucg}{{\cal U}_{h,r_{cg}}({\frak g})}
\nc{\Ua}{{\cal U}_{h,a}({\frak g})}
\nc{\Zg}{{\cal Z}({\frak g})}
\nc{\Zk}{{\cal Z}({\frak k})}
\nc{\Sh}{{\cal S}({\frak h})}
\nc{\Uh}{{\cal U}({\frak h})}
\nc{\Ut}{{\cal U}({\frak t})}
\nc{\Uk}{{\cal U}({\frak k})}
\nc{\cA}{\cal A}
\nc{\cS}{\cal S}
\nc{\cC}{\cal C}
\nc{\cZ}{\cal Z}
\nc{\bZ}{\overline{\cal Z}}
\nc{\bZp}{\overline{\cal Z}_{\bullet}}
\nc{\cZp}{{\cal Z}_{\bullet}}
\nc{\wZ}{\widetilde{\cal Z}}
\nc{\cG}{\cal G}
\nc{\cL}{\cal L}
\nc{\cR}{\cal R}
\nc{\cU}{\cal U}
\nc{\cH}{\cal H}
\nc{\wH}{\widetilde{\cal H}}
\nc{\bH}{\overline{\cal H}}
\nc{\bHp}{\overline{\cal H}_{\bullet}}
\nc{\cK}{\cal K}
\nc{\cF}{\cal F}
\nc{\cE}{\cal E}
\nc{\fg}{\frak g}
\nc{\CO}{\cal O}
\nc{\fn}{\frak n}
\nc{\fm}{\frak m}
\nc{\fsl}{\frak sl}
\nc{\fh}{\frak h}
\nc{\fb}{\frak b}
\nc{\ft}{\frak t}
\nc{\fk}{\frak k}
\nc{\fp}{\frak p}
\nc{\fI}{\frak I}
\nc{\C}{\Bbb C}
\nc{\Q}{\Bbb Q}
\nc{\Z}{\Bbb Z}
\nc{\N}{\Bbb N}

\nc{\dirlim}{\underset{\rightarrow}{\lim}\,} 
\nc{\nen}{\newenvironment}
\nc{\ol}{\overline}
\nc{\ul}{\underline}
\nc{\ra}{\rightarrow}
\nc{\lra}{\longrightarrow}
\nc{\Lra}{\Longrightarrow}
\nc{\Lla}{\Longleftarrow}
\nc{\Llra}{\Longleftrightarrow}
\nc{\hra}{\hookrightarrow}
\nc{\iso}{\overset{\sim}{\lra}}
\nc{\ssubset}{\underset{\not=}{\subset}}

\nc{\Thm}[1]{Theorem~\ref{#1}}
\nc{\Thme}[1]{Th\'eor\`eme~\ref{#1}}
\nc{\Prop}[1]{Proposition~\ref{#1}}
\nc{\Lem}[1]{Lemma~\ref{#1}}
\nc{\Cor}[1]{Corollary~\ref{#1}}
\nc{\Conj}[1]{Conjecture~\ref{#1}}
\nc{\Claim}[1]{Claim~\ref{#1}}
\nc{\Defn}[1]{Definition~\ref{#1}}
\nc{\Exa}[1]{Example~\ref{#1}}
\nc{\Rem}[1]{Remark~\ref{#1}}
\nc{\Note}[1]{Note~\ref{#1}}
\nc{\Quest}[1]{Question~\ref{#1}}
\nc{\Hyp}[1]{Hypoth\`ese~\ref{#1}}

\nen{thm}[1]{\label{#1}{\bf Theorem.\ } \em}{}
\nen{thme}[1]{\label{#1}{\bf Th\'eor\`eme.\ } \em}{}
\nen{prop}[1]{\label{#1}{\bf Proposition.\ } \em}{}
\nen{lem}[1]{\label{#1}{\bf Lemma.\ } \em}{}
\nen{cor}[1]{\label{#1}{\bf Corollary.\ } \em}{}
\nen{conj}[1]{\label{#1}{\bf Conjecture.\ } \em}{}
\nen{claim}[1]{\label{#1}{\bf Claim.\ } \em}{}

\nen{defn}[1]{\label{#1}{\bf Definition.\ } }{}
\nen{exa}[1]{\label{#1}{\bf Example.\ } }{}

\nen{rem}[1]{\label{#1}{\em Remark.\ } }{}
\nen{note}[1]{\label{#1}{\em Note.\ } }{}
\nen{exer}[1]{\label{#1}{\em Exercise.\ } }{}
\nen{sket}[1]{\label{#1}{\em Sketch of proof.\ } }{}
\nen{quest}[1]{\label{#1}{\bf Question.\ } \em}{}
\nen{hyp}[1]{\label{#1}{\bf Hypoth\`ese.\ } \em}{}

\newenvironment{guillemet}{$\scriptstyle <<$\kern.1em}
{\kern.1em$\scriptstyle >>$}

\begin{document}

\title[]{The Zhang transformation and 
$\cU_q\bigl(\osp(1,2l)\bigr)$-Verma modules annihilators}

\author[]{Emmanuel Lanzmann}
\address{ 
Dept. of Theoretical Mathematics, The Weizmann Institute of Science, 
Rehovot 76100, Israel,
{\tt email: lanzmann@@wisdom.weizmann.ac.il} 
}

\thanks{
The author was partially supported by the EC TMR network Algebraic
Lie Representations 
Grant No. ERB FMRX-CT97-0100 and
Minerva grant 8337}

\begin{abstract} 
In~\cite{zh}, R.~B.~Zhang found a way to link certain formal deformations 
of the Lie algebra $\oo(2l+1)$ and the Lie superalgebra $\osp(1,2l)$. 
The aim of this article is to reformulate the Zhang
transformation 
in the context of the quantum enveloping algebras 
{\em \`a la}  Drinfeld-Jimbo and to 
show how this construction can explain
 the main theorem of~\cite{gl2}:
the annihilator of a Verma module over
the Lie superalgebra $\osp(1,2l)$ is 
generated by its intersection with the centralizer of the even part
of the enveloping algbra.
\end{abstract}
\maketitle

\section{Introdution}

A well known theorem of Duflo claims that the 
 annihilator of a Verma module over a complex semi-simple Lie algebra 
is 
generated by its intersection with the centre of the enveloping algbra.
In~\cite{gl2} we show that in order for this theorem to hold in the case
of the Lie superalgebra $\osp(1,2l)$ one has to replace 
the centre by the centralizer of the even part
of the enveloping algbra.
The purpose of this article is to show how quantum groups can illucidate this
phenomemon.

Let $\fk$, $\fg$ be respectively 
the complex simple Lie algebra $\oo(2l+1)$ and the complex superalgebra
$\osp(1,2l)$.  From many point of
views, the algebras $\fg$ and $\fk$ are very similar objects. For instance,
identifying properly 
Cartan subalgebras of $\fk$ and $\fg$, 
the root systems $\Delta_{\fk}$, $\Delta_{\fg}$ are contained 
one into the other, and the set of irreductible roots of $\Delta_{\fg}$ 
is $\Delta_{\fk}$. 
Moreover, given a simple finite dimensional $\fg$-module, 
the corresponding simple $\fk$-module of the same highest weight
is also finite dimensional 
and has the same formal character (and even the same crystal graph).
 Nevertheless, there is no obvious
direct way to link the algebras $\fg$, $\fk$. 
To bridge the gap, one has to go through the quantum level:
in his article published in 1992, R.~B.~Zhang (see~\cite{zh}, 3) found
a recipe to pass from a certain formal deformation of  $\cU(\fk)$ to
a formal deformation of
$\cU(\fg)$.

In  this article we present a  reformulation of  the 
Zhang transformation in the more algebraic  context of the 
quantizations {\it \`a la} Drinfeld-Jimbo. 
The idea is to start with the Drinfeld-Jimbo
quantum enveloping algebra 
$\cU:=\cU_{\sqrt q}(\oo(2l+1))$ and to extend the torus by
the finite group $\Gamma:=\Delta_{\fg}/ 2\Delta_{\fg}$. In other 
words, we introduce the semi-direct product 
$\wU:=\cU\rtimes k\Gamma$  
where $\Gamma$ acts on $\cU$ in an obvious
manner. Twisting the generators of $\cU$ by elements of $\Gamma$, 
we build a subalgebra $\bU$ of $\wU$ isomorphic to the quantum
enveloping algebra $\cU_{-q}(\fg)$,
and such that $\wU\simeq \bU\rtimes k\Gamma$. This 
construction provides an involution of the vector space $\wU$
mapping $\cU$
onto $\bU$. We call it the  Zhang transformation. 

A first obvious consequence of this construction is that the 
algebras $\cU$ and $\bU$ have the ``same''
 finite dimensional modules. To be more precise,
the quantum simple spinorial $\cU$-modules, viewed as $\bU$-modules (via 
$\wU$), 
are not deformation of $\fg$-modules (even up to 
 tensorization by one dimensional modules). Moreover, roughly speaking,
given a  simple finite dimensional $\cU$-modules which is not of spinorial
type, the specialization $q\rightarrow 1$ provides a simple $\fk$-module and 
the specialization $q\rightarrow -1$ a simple $\fg$-module. These 
classical $\fk$ and $\fg$-modules
  have therefore the same formal character. We believe that the characters
of their Demazure modules are also equal, and this construction might be an
interesting approach to this problem.

Another consequence, which is the crucial observation for our purpose, is that 
the Zhang transformation maps the centre of $\cZ(\cU)$
to $\cA(\bU)$, the commutant of the even part of the $\bU$. 
 As in the classical case, $\cA(\bU)$ turns out to be 
 the direct sum of the centre $\cZ(\bU)$ of $\bU$ and of the anticentre
$\cZp(\bU)$, the subspace of elements  which commute with even elements 
and anticommute with odd elements.  
The subspace $\cZp(\bU)$ is a cyclic module over the centre 
$\cZ(\bU)$. We construct a generator of this module which is 
a quantization of the element $T$ introduced in 4.4.1~\cite{gl2}. 
The set  $\cZp(\bU)$ has 
another interpretation: this is  the set of invariant elements under
a certain twisted superadjoint action.  This twisted action is the 
quantum analogue of the ``non-standard'' 
action considered in the classical case by Arnaudon, Bauer, Frappat 
(see~\cite{abf}, 2). 
More generally, we show that the locally finite part $\F(\cU)$ of
$\cU$ for the adjoint action is mapped to the direct sum
$\bF(\bU)\oplus \bFp(\bU)$ of the locally finite parts of $\bU$ for
the superadjoint action  and for its twisted version. 
This allows us to deduce from
the work of Joseph and Leztzer (see~\cite{jl})
 that  the annihilator of a $\bU$-Verma 
module in  
$\bF(\bU)\oplus \bFp(\bU)$ is generated by its intersection with
$\cA(\bU)$ (theorem~\ref{qthmduflo}).
We also prove a quantum analogue of  theorem 7.1~\cite{gl}. 

The article is organized as follows. In section~\ref{qsecwU} we present
the algebra $\wU$, the main object of this article, and we define
the Zhang transformation. In section~\ref{qtwisad}
we study the locally finite parts of $\wU$
for different actions and analyse how these actions are transformed
by the Zhang transformation. We define in~\ref{secalgstrucbU} the
subalgebra $\bU$. We show that $\bU$ is isomorphic to $\cU_{-q}(\fg)$ 
and we deduce from  section~\ref{qtwisad} algebraic structure 
theorems for $\bU$. Section~\ref{qsecdufloth} is devoted to a proof
of the annihilation theorem for $\bF(\bU)\oplus \bFp(\bU)$.
In section~\ref{qseccom} we show that 
$\cA(\bU)$ coincides with the centralizer of the even part of $\bU$. 

{\em Acknowledgement}. I am grateful to M.~Gorelik for reading  earlier
versions of this paper and making numerous important remarks. I would like
to thank T.~Joseph for his support and his comments. I wish to express
my gratitude to P.~Cartier, M.~Duflo, T.~ Levasseur and G.~Perets for
helpful discussions.

\vfill\break

\section{Background}
\subsection{Notations}
\label{notation4}

Let ${\frak k}$ be the complex simple Lie algebra $\oo(2l+1)$, $l\geq 1$.
Fix a Cartan subalgebra
${\frak h}$ of $\fk$ and denote by $\Delta_{\fk}$ the root system of $\fk$.
We fix a basis of simple roots $\pi$ of $\Delta_{\fk}$.
 Denote by $W$ the Weyl group of $\Delta_{\fk}$, and set 
$\rho:=\displaystyle\mathop{\sum}_{\alpha\in
\Delta_{\fk}^+}\alpha$.
Denote by $(-,-)$ the non-degenerate bilinear form on ${\frak h}^*$ coming 
from the
restriction of the Killing form of ${\frak g}_0$
to ${\frak h}$.
For any $\lambda,\mu\in 
{\frak h}^*,\
(\mu,\mu)\not=0$ one defines 
$\langle\lambda,\mu\rangle:=2(\lambda,\mu)/ (\mu,\mu).$
One has the following useful realization of $\Delta_{\fk}$. 
Identify ${\frak h}^*$ with ${\Bbb C}^l$ and consider $(-,-)$ as 
an inner product on ${\Bbb C}^l$. Then there exists an orthonormal basis 
$\{\beta_1,\ldots,\beta_l\}$ such that
$$\pi=\{\beta_1-\beta_2,\ldots,\beta_{l-1}-\beta_l,\beta_l\},\ \  
\Delta_{\fk}=
\{\pm\beta_i\pm\beta_j, 1\leq i<j\leq l,\ \pm\beta_i,\ 1\leq i\leq l\}$$
and the Weyl group $W$ is just the group of the signed permutations of the 
$\beta_i$. We set $\alpha_i:=\beta_{i}-\beta_{i+1}$ with
the convention $\beta_{l+1}=:0$. Let also $w_i$,  $1\leq i\leq l$, be
the elements $w_i:=\beta_1+\ldots +\beta_i$. With these notations,
the set of weights $P_{\fk}(\pi)$ of $\Delta_{\fk}$ is 
$$ P_{\fk}(\pi):=\Z w_1\oplus\ldots \Z w_{l-1}\oplus \Z (w_l/ 2).$$
We consider also the set of dominant weights $ 
P_{\fk}^+(\pi)$.

\subsubsection{} Let $\fg$ be the complex Lie superalgebra $\osp(1,2l)$. 
Denote
by ${\frak g}_0$ the even part of $\fg$ and by ${\frak g}_1$ its odd part.
We identify $\fh$ with a Cartan subalgebra of $\fg$ in such a way
that $\pi$ is also a basis of simple roots of the root system $\Delta_{\fg}$
of $\fg$ with respect to $\fh$. 
The sets of even and odd roots of $\fg$  equal
respectively $\{\pm\beta_i\pm\beta_j, 1\leq i<j\leq l,\ \pm 2\beta_i,
1\leq i\leq l\}$ and $\{\pm \beta_i,\ 1\leq i\leq l\}$. The set of 
weights $P_{\fg}(\pi)$ of $\Delta_{\fg}$ is 
$$P_{\fg}(\pi):=\Z w_1\oplus\ldots \Z w_{l-1}\oplus \Z w_l
=\Z \alpha_1\oplus\ldots \Z \alpha_{l-1}\oplus \Z \alpha_l.$$
In a 
standard manner we define the set  of dominant
weights $P_{\fg}^+(\pi)$.

\subsubsection{}\label{qcoef}
Let $q$ be an indeterminate and  $k=\C(\sqrt{q})$. Let $\nu\in k$, be
such that 
$\nu\not= 0,\pm 1$. For all $n\in \N$, we set  
$[n]_{\nu}:=\displaystyle{\nu^n-\nu^{-n}\over 
\nu-\nu^{-1}}$ and $[n]_{\nu}!:=[n]_{\nu}\times[n-1]_{\nu}
\times\ldots\times [1]_{\nu}$ with the convention $[0]_{\nu}=1$.
If $1\leq m\leq n\in \N$ we set $\binome{\nu}nm:=
\displaystyle{[n]_{\nu}!\over [m]_{\nu}![n-m]_{\nu}!}$.

\subsection{The algebra $\cU$}\label{qsecdefU}

 Let  $\cU$ be the algebra over the field $k$ generated by
 the elements $E_i$, $F_i$,  $1\leq i\leq l$,  
$K_{\mu}$,   $\mu\in P_{\fg}(\pi)$ under the relations 
\begin{equation}\label{defwueq1}
 K_0=1,\ K_{\lambda}K_{\mu}=K_{\lambda+\mu}
\end{equation}
\begin{equation}\label{defwueq3}
K_{\lambda}E_jK_{\lambda}^{-1}=q^{(\lambda,\alpha_j)}E_j,\ 
K_{\lambda}F_jK_{\lambda}^{-1}=
q^{-(\lambda,\alpha_j)}F_j
\end{equation}
\begin{equation}\label{reluhef}\label{defwueq4}
E_iF_j-F_jE_j=\delta_{ij}\displaystyle{ K_{\alpha_i}-K_{\alpha_i}^{-1}
\over q-q^{-1}}
\end{equation}
\begin{equation}\label{reluhsere}
\displaystyle\sum_{k=0}^{1-\langle\alpha_j,\alpha_i\rangle}{(-1)}^k
\binome {q_i}{1-\langle\alpha_j,\alpha_i\rangle}k
E_{i}^{1-\langle\alpha_j,\alpha_i\rangle-k}
E_{j}E_{i}^{k}=0
\end{equation}
\begin{equation}\label{reluhserf}
\displaystyle\sum_{k=0}^{1-\langle\alpha_j,\alpha_i\rangle}{(-1)}^k
\binome {q_i}{1-\langle\alpha_j,\alpha_i\rangle}k
F_{i}^{1-\langle\alpha_j,\alpha_i\rangle-k}
F_{j}F_{i}^{k}=0
\end{equation}
where $q_i:={q}^{(\alpha_i,\alpha_i)\over 2}$ and 
$\langle\alpha_j,\alpha_i\rangle:= 2(\alpha_j,\alpha_i)/
(\alpha_i,\alpha_i)$.
Relations (\ref{reluhsere}),
 (\ref{reluhserf}) are called the 
quantum  Serre relations. 

\subsubsection{}\label{qsubsecU}
Replacing $E_l$ by 
$\left({\displaystyle{q-q^{-1}}\over \sqrt{q}-\sqrt{q}^{-1}}\right)
E_l$ in the above equations, we see that $\cU$
 is just the Drinfeld-Jimbo quantum enveloping algebra  
$\cU_{\sqrt q}(\fk)$.

\subsubsection{}\label{qnottorus}
Let $\cU^+$ (resp. $\cU^-$) be the subalgebra of  $\cU$
generated by the $E_i$ (resp. $F_i$).
We also denote by $\cU^o$ the subalgebra generated
by the $K_{\lambda}$. If $T$ stands for the multiplicative group
$\{K_{\mu},\ \mu\in P_{\fg}(\pi)\}$, then the
group algebra of $T$ identifies with $\cU^o$. 
One has the triangular decomposition
$\wU\simeq \cU^-\otimes \cU^o \otimes \cU^+$,
this isomorphism of vector spaces being given by multiplication.

\subsubsection{}\label{qHopf}
The algebra $\cU$ is a  $k$-Hopf algebra with counit
$\varepsilon$, coproduct $\Delta$ and antipode $S$  defined by 
\begin{equation*}
\varepsilon(K_{\lambda})=1,\ 
\varepsilon(E_i)=\varepsilon(F_i)=0
\end{equation*}
\begin{equation*}
\begin{array}{c}\Delta(K_{\lambda})=K_{\lambda}\otimes K_{\lambda},\\
 \Delta(E_i)=E_i\otimes 1+K_{\alpha_i}\otimes E_i,\
\Delta(F_i)=F_i\otimes K_{\alpha_i}^{-1}+1\otimes F_i\\
\end{array}
\end{equation*}
\begin{equation*}
S(K_{\lambda})=K_{\lambda}^{-1},\
S(E_i)=-K_{\alpha_i}^{-1}E_i,\ 
S(F_i)=-F_iK_{\alpha_i}
\end{equation*}

\subsubsection{Representations}\label{repU}
Let
$\widehat{T}$ be
the group of characters of $T$ with values in $k$. 
Given a character $\Lambda\in
\widehat{T}$,
we say that $\Lambda$ is linear if there exists $\lambda\in P_{\fk}(\pi)$
such that $\Lambda(K_{\mu})=q^{(\lambda,\mu)}\ \forall \mu\in 
P_{\fg}(\pi)$. In that case we write
$\Lambda:=q^{\lambda}$.   

Let $M$ be a $T$-module. An element $m\in M$
is said to be an  element of weight $\Lambda\in \widehat{T}$ if
$K_{\mu}m=\Lambda(K_{\mu})m\ \forall \mu\in P_{\fg}(\pi)$. We also
call $m$ a $T$-weight element. We denote by $M_{\Lambda}$ 
the subspace of elements of weight $\Lambda$. 

For any $\Lambda\in \widehat{T}$, we denote by $M(\Lambda)$ the 
$\cU$-Verma module of highest weight $\Lambda$, and by $V(\Lambda)$
its unique simple quotient. We recall now basic properties of
the representation theory of $\cU$ 
 (see for instance~\cite{jq}, 4.4.9). 
 If $\Lambda=q^{\lambda}$ is linear,  
and if there exists $\alpha\in \Delta^+_{\fk}$ such that 
$\langle \lambda+\rho,\alpha\rangle\in \N^+$ then 
$M(q^{s_{\alpha}.\lambda})$ is 
a submodule of $M(q^{\lambda})$. 
The simple module $V(\Lambda)$
is finite dimensional if and only if $\Lambda=\phi q^{\lambda}$ where
$\lambda\in P_{\fk}^+(\pi)$ and 
$\phi\in  \widehat{T}$ is such that $\phi(K_{\mu})=\pm 1$ for all 
$\mu\in P_{\fg}(\pi)$.  
Any  finite dimensional $\cU$-module 
is completely reducible. 

\subsubsection{}\label{weightU}
 The group $T$ acts on $\cU$ by inner automorphisms. 
Thus we can speak of weight elements in $\cU$. By a slight abuse of 
notation we shall say that $u\in \cU$ is of weight 
$\lambda\in P_{\fk}(\pi)$ if it is actually of weight $q^{\lambda}$.

\section{The algebra $\wU$}\label{qsecwU}
In this section we  introduce the main object of this article.

\subsection{Definition and basic properties}
 \subsubsection{Definition}\label{qdefwU}
 Let $\Gamma$ be the multiplicative group
$\{\xi_{\mu}, \mu\in P_{\fg}(\pi)\}\simeq P_{\fg}(\pi)\slash 2P_{\fg}(\pi)$
and $k\Gamma$ its group algebra. The group $\Gamma$ acts on 
$\cU$ in the following natural way:
$$\xi_{\mu}. u ={(-1)}^{(\mu,\lambda)}u,\ \ \forall u\in\cU \mbox{ of
weight $\lambda$},\ \forall \lambda,\mu\in P_{\fg}(\pi).$$ 
Using this action of $\Gamma$ we introduce the
algebra 
$$\wU:=\cU\rtimes k \Gamma.$$
Throughout
this article, we shall use
the shortened notation $u\xi_{\mu}$ for an
 element $u\otimes \xi_{\mu}\in \wU$. 

\subsubsection{The Hopf structure}\label{HopfwU} Recall that the
 algebra $k\Gamma$ is a Hopf algebra for the 
counit $\varepsilon$, coproduct $\Delta$ and antipode $S$
defined by $$\varepsilon(\xi_{\mu})=1,\ 
\Delta(\xi_{\mu})=\xi_{\mu}\otimes \xi_{\mu},\ S(\xi_{\mu})=\xi_{\mu}.$$
Since $\Gamma$ acts on $\cU$ by Hopf algebra automorphisms,
$\wU=\cU\rtimes k \Gamma$ is a Hopf algebra for the obvious 
coalgebra 
structure (resp. antipode) on 
tensor product.

\subsubsection{A $\Z_2$-gradation}
\label{qgrad} The algebra $\wU$ is also endowed with the following 
$\Z_2$-gradation. 
Introduce $\xi:=\xi_{w_l}$, and define  
$$\ \forall \overline{i}\in\Z_2,\ \ {\wU}_{|_i}=\{x\in \wU,\
\xi x\xi={(-1)}^{i}x\}.$$
An element $x\in {\wU}_{|_0}\cup {\wU}_{|_1} $ is called $\Z_2$-homogenous,
 and we write $|x|=i$ if $x\in {\wU}_{|_i}$, $x\not=0$.

\subsubsection{The Harish-Chandra projection}
\label{defhcq} Recall the triangular decomposition given 
in~\ref{qnottorus}. Introduce the subalgebra $\wU^{o}:=\cU^o\otimes
k\Gamma$. One has the triangular decomposition
$$\wU\simeq \cU^-\otimes \wU^o\otimes \cU^+.$$ 
Let $\cU^{++}$ (resp. $\cU^{--}$) be the augmentation ideal of 
$\cU^+$ (resp. $\cU^-$). The triangular decomposition of $\wU$
implies $\wU=(\cU^{--}\wU+\wU\cU^{++})\oplus
\wU^0$ 
which allows to define a Harish-Chandra projection 
$\Upsilon: \wU\longrightarrow  \wU^o$ 
with respect to this decomposition.

\subsection{Representations of $\wU$}
\subsubsection{Generalities}\label{qrepgen}
 Let  
$\widehat{\Gamma}$ be the group of characters of $\Gamma$, identified
with the set of group morphisms  $P_{\fg}(\pi)\slash
2P_{\fg}(\pi)\longrightarrow \{1,-1\}$. Any $\lambda\in 
P_{\fg}(\pi)$ defines a character 
$(-1)^{\lambda}\in \widehat{\Gamma}$ by the formula: 
$(-1)^{\lambda}(\mu):={(-1)}^{(\lambda,\mu)}.$
 Observe
that $\Gamma$ embeds  in $\widehat{T}$ as the set 
$\{\Lambda\in \widehat{T}\mbox{ s.t. }\ \forall \mu\in P_{\fg}(\pi),\ 
\Lambda(K_{\mu})=\pm 1\}$.

Recall~\ref{repU}. Let $M$ be a $T\times \Gamma$-module. We say that an element
$m\in M$ is a $T\times \Gamma$-weight element of weight $(\Lambda,\theta)
\in\widehat{T}\times \widehat{\Gamma}$ if 
$\xi_{\mu'}K_{\mu}m=\Lambda(K_{\mu})\theta(\mu')m\
 \forall \mu,\mu'\in P_{\fg}(\pi)$. 

Take $\Lambda\in\widehat{T}$
and $\theta\in\widehat{\Gamma}$. There is an obvious way to endow
$M(\Lambda)$ with a structure of a $\wU$-module. Define for any
$x\in  M(\Lambda)$ of weight $q^{-\nu}\Lambda$, $\nu\in P_{\fg}(\pi)$,
$$\xi_{\mu}x:={(-1)}^{(\mu,\nu)}\theta(\mu)x\ \ \forall \mu\in P_{\fg}(\pi).$$
We call this $\wU$-module a $\wU$-Verma module and we denote
it by $M(\Lambda,\theta)$. By definition, 
$M(\Lambda,\theta)$ and $M(\Lambda)$ have the same submodules. 
Let $\Lambda=q^{\lambda}$ be linear. Assume that   
there exists $\alpha\in \Delta^+_{\fk}$ such that 
$\langle \lambda+\rho,\alpha\rangle\in \N^+$
and define $\theta':={(-1)}^{s_{\alpha}.\lambda-\lambda}\theta$.
 Then 
$M(q^{s_{\alpha}.\lambda},\theta')$ is 
a $\wU$-submodule of $M(q^{\lambda},\theta)$.
 The $\wU$-module
$M(\Lambda,\theta)$ has a unique simple quotient, $V(\Lambda,\theta)$.
As a $\cU$-module, $V(\Lambda,\theta)\simeq V(\Lambda)$.

\subsubsection{} \label{propcomplredwU}
\begin{lem}{} All finite dimensional  $\wU$-modules are completely
  reducible. Moreover, any simple finite dimensional $\wU$-module 
$M$ is isomorphic to 
a $V(q^{\lambda}\phi,\theta)$, $\lambda\in  P^+_{\fk}(\pi)$,
$\phi,\theta\in \widehat{\Gamma}$.
\end{lem}
\begin{pf} Let $M^{{\cU}^+}$ be the subspace of $M$ of invariant vectors by
$\cU^+$. This subspace is stable by the action of the commutative 
algebra $\wU^o$. Since $\Gamma$ is finite and $M$ is a $\cU$-module of 
finite dimension, $\wU^o$
 acts diagonally on $M$ and hence on $M^{{\cU}^+}$. 
Let $\{v_1,\ldots,v_r\}$ be basis
of $M^{{\cU}^+}$ composed of $T\times \Gamma$-weight vectors. The
representation
theory of $\cU$ (see~\cite{jan} chap. 5) asserts that on the one hand 
$M_i:=\cU v_i=\wU v_i$ is a simple $\cU$-module, and so a
$\wU$-module 
of the form  $V(q^{\lambda_i}\phi_i,\theta_i)$ with  
$\lambda_i\in P^+_{\fk}(\pi)$,  $\theta,\phi_i\in \widehat{\Gamma}$,
and
on the other hand that  $M=\oplus M_i$.
\end{pf}

\subsubsection{The $\Z_2$-gradations}\label{qrepgengradM}
Both $M(\Lambda,\theta)$ and $V(\Lambda,\theta)$ inherits a 
$\Z_2$-gradation. Let $v_{\Lambda}$ be the highest weight vector of
 $M(\Lambda,\theta)$. They are two natural 
 $\Z_2$-gradations on $M(\Lambda,\theta)$. Fix  ${j}\in \Z_2$ and
  define  ${M(\Lambda,\theta)}_{|_{ i}}=
\wU_{|_{ {i}+{j}}}v_{\lambda}\ \forall
{i}\in \Z_2$. The gradations on  $V(\Lambda,\theta)$ are defined
similarly using the highest weight vector of
$V(\Lambda,\theta)$.

\subsection{Three gradations and the Zhang transformation}
\label{qsecUUs}

\subsubsection{The gradation by the weights}
\label{notweight}  The considerations of~\ref{weightU} extend to $\wU$. 
We shall denote by $\nu(x)$ the weight of 
a weight element $x\in\wU$ and by $\wU_{\nu}$ the subspace of elements 
of weight $\nu$. The algebra $\wU$ is graded by its weight subspace:
$$\wU:=\bigoplus_{\nu\in P_{\fg}(\pi)} \wU_{\nu}.$$

\subsubsection{The $\mu$-gradation and the Zhang transformation}
\label{defzhtransfor}
Define the $P_{\fg}(\pi)/2P_{\fg}(\pi)$-gradation 
$$\wU:=\bigoplus_{\mu\in P_{\fg}(\pi)/2P_{\fg}(\pi)}{}^{\mu}\wU$$
for which $\xi_{\lambda}\in {}^{0}\wU$, $K_{\lambda}\in {}^{\lambda}\wU$,
 $E_i\in {}^{\beta_{i+1}}\wU$, and
$F_i\in{}^{\beta_{i}}\wU$ for all $\lambda\in P_{\fg}(\pi)$, $1\leq i\leq l$.
If $x\in {}^{\mu}\wU$, we set $\mu(x):=\mu$.   
We call  Zhang transformation the involution of vector space
$\Psi:\wU\longrightarrow \wU$ such that 
$$\Psi(x):=\xi_{\mu} x\ \ \forall
x\in {}^{\mu}\cU,\ \forall \mu\in P_{\fg}(\pi)/2P_{\fg}(\pi).
$$ We shall show in~\ref{qalgstUs} that $\Psi(\cU)$ is 
a subalgebra isomorphic to $\cU_{-q}(\fg)$ (see~\cite{mz} for
the definition of this algebra).
For any   homogenous elements $a,b$ for the respective gradations 
$(\wU_{\nu})$ and $({}^{\mu}\wU)$,
\begin{equation}\label{qcomPsi}
\Psi(ab)={(-1)}^{(\nu(a),\mu(b))}\Psi(a)\Psi(b).
\end{equation} 
\subsubsection{The $\delta$-gradation}\label{qdefgraddelta}
We introduce another $P_{\fg}(\pi)/2 P_{\fg}(\pi)$-gradation on $\wU$ 
(compare this gradation with the filtration defined in 5.3.1~\cite{jq}) 
$$\wU=\bigoplus_{\delta\in P_{\fg}(\pi)/2 P_{\fg}(\pi) } \wU^{\delta}$$
 for which $\xi_{\mu},E_i\in \wU^0$, 
$K_{\mu}\in \wU^{\mu}$, $F_i\in \wU^{\alpha_i}$ for
all $\mu\in P_{\fg}(\pi)$, $1\leq i\leq l$. A glance at
the defining relations of $\wU$ ensures that this does define a gradation
on $\wU$. 
If $x\in \wU^{\delta}$, we shall use the notation $\delta(x):=\delta$.

\subsubsection{Compatibilities between the gradations}
The gradation $(\wU^{\delta})$
is invariant under the action of $T$ (see~\ref{notweight}). Thus the
$\wU^{\delta}$ are direct  sums of there weight subspaces.  
If $\delta,\nu\in P_{\fg}(\pi)$, we set $\wU_{\nu}^{\delta}:=\wU^{\delta}
\cap \wU_{\nu}$. One has the bigradation on $\wU$
\begin{equation}
\label{eqbigrad}
\wU=\bigoplus_{\nu\in P_{\fg}(\pi)\atop
\delta\in P_{\fg}(\pi)/2 P_{\fg}(\pi)}\wU_{\nu}^{\delta}.
\end{equation}
The relation between the gradations $(\wU_{\nu})$, $({}^{\mu}
\wU)$, $(\wU^{\delta})$ reads as follows. The gradation $({}^{\mu}
\wU)$ is also $T$-invariant, and hence induces a bigradation 
$\wU=\oplus \bigl({}^{\mu}\wU\cap \wU_{\nu})$. Then, this bigradation 
coincides
with the bigradation (\ref{eqbigrad}). To be more precise, one has 
\begin{equation}\label{qcompagrad}
 \wU_{\nu}^{\delta}= {}^{\delta+ \eta(\nu)}\wU\cap \wU_{\nu}
\end{equation} 
where $\eta: P_{\fg}(\pi)\rightarrow  P_{\fg}(\pi)/2P_{\fg}(\pi)$ is
the map
defined by $$\eta( \sum n_i\alpha_i):=\sum n_i\beta_{i+1}.$$
Indeed, it is enough to check (\ref{qcompagrad}) on the generators. 
One has $\delta(E_i)+\eta(\nu(E_i))=\beta_{i+1}=\mu(E_i)$,
$\delta(F_i)+\eta(\nu(F_i))=\alpha_i+\beta_{i+1}=\beta_i=\mu(F_i)$,
$\delta(K_{\mu})+\eta(\nu(K_{\mu}))=\mu=\mu(K_{\mu})$,
$\delta(\xi_{\mu})=\nu(\xi_{\mu})=\mu(\xi_{\mu})=0$.

\subsubsection{}\label{qparity} Recall the $\Z_2$-gradation defined 
in~\ref{qgrad}. If $\nu\in P_{\fg}(\pi)$, we set $|\nu|:=(\nu,w_l)\pmod 2$.
Observe that for all $\nu,\nu'\in P_{\fg}(\pi)$, the following 
identity holds in $\Z_2$:
\begin{equation}\label{qeqparity}
(\nu,\eta(\nu'))+(\eta(\nu),\nu')+(\nu,\nu')=|\nu||\nu'|
\end{equation}
Since both sides of (\ref{qeqparity}) are bilinear in $\nu,\nu'$, 
the identity (\ref{qeqparity}) reduces to the case where $\nu=\alpha_i$,
$\nu'=\alpha_j$. In that case, the left hand side of (\ref{qeqparity})
is equal (in $\Z_2$) to 
$$(\alpha_i,\beta_{j+1})+(\beta_{i+1},\alpha_j)+(\alpha_i,\alpha_j)
= (\beta_i,\beta_j)+(\beta_{i+1},\beta_{j+1}) =\left\{\vcenter{
\hbox{$0$ if $(i,j)\not= (l,l)$}
\hbox{$1$ if $i=j=l$}}\right.= |\alpha_i||\alpha_j|$$

\subsubsection{}\label{qcomHopfbigrad}
Recall the definition of $\wU_{\nu}^{\delta}$ 
(see~\ref{qdefgraddelta}). One has 
\begin{equation}\label{qcomcopbigrad}
\begin{array}{l} \Delta (\wU^{\delta}_{\nu})\subset
\displaystyle\bigoplus_{\nu_1+\nu_2=\nu}
\wU^{\delta+\nu_2}_{\nu_1}\otimes \wU^{\delta}_{\nu_2}\\
S(\wU^{\delta}_{\nu})\subset \wU^{\delta+\nu}_{\nu}
\end{array}
\end{equation}
Indeed, according to~\ref{qHopf},~\ref{HopfwU} these inclusions
are satisfied for the generators of $\wU$.

\subsection{A Hopf superalgebra structure on $\wU$}\label{sHopfwU}
The algebra $\Psi(\cU)$ is not a Hopf subalgebra of 
$\wU$ for the coproduct and antipode defined in~\ref{qHopf},~\ref{HopfwU}.
In this subsection, we   endow $\wU$ with a structure of 
Hopf superalgebra for which $\Psi(\cU)$ is  a Hopf subalgebra.

Recall definitions of $\Psi$ (see~\ref{defzhtransfor}) and 
$\wU_{\nu}^{\delta}$ 
(see~\ref{qdefgraddelta}). Define for all homogenous elements $a\in \wU$
for the bigradation $\wU_{\nu}^{\delta}$:
\begin{equation}\label{qeqdefsHopf}
\begin{array}{lcl} 
\overline{\Delta} \Psi(a)&:=& {(-1)}^{(\nu(a_1),\nu(a_2)+\eta(\nu(a_2))}
\Psi(a_1)\otimes \Psi(a_2)\\
\overline{S}\Psi(a)&:=& {(-1)}^{(\nu(a),\delta(a))}\Psi(S(a))
\end{array}
\end{equation}
with  Sweedler notation $\Delta(a)=a_1\otimes a_2$. 

\subsubsection{}\label{qdefsgen}
Intoduce $e_i:=\Psi(E_i)$, $f_i=\Psi(F_i)$, $k_{\mu}:=\Psi(K_{\mu})$. 
On the generators the definitions (\ref{qeqdefsHopf}) give:
\begin{equation*}
\begin{array}{c}\overline{\Delta}(\xi_{\mu})=\xi_{\mu}\otimes \xi_{\mu},
\ \overline{\Delta}(k_{\lambda})=
k_{\lambda}\otimes k_{\lambda},\\
 \overline{\Delta}(e_i)=e_i\otimes 1+k_{\alpha_i}\otimes e_i,\
\overline{\Delta}(f_i)=f_i\otimes k_{\alpha_i}^{-1}+1\otimes f_i\\
\end{array}
\end{equation*}
\begin{equation*} \overline{S}(\xi_{\mu})=\xi_{\mu},\ 
\overline{S}(k_{\lambda})=k_{\lambda}^{-1},\
\overline{S}(e_i)=-k_{\alpha_i}^{-1}e_i,\ 
\overline{S}(f_i)=-f_ik_{\alpha_i}
\end{equation*}
\subsubsection{}\begin{lem}{} $\wU$ endowed with $(\overline{\Delta},
\overline{S}, \varepsilon)$ is a Hopf superalgebra.
\end{lem}
\begin{pf}
 Retain the definition of the $\Z_2$-gradation on $\wU$ 
(see~\ref{qgrad}). 
Let us prove that for any $a,b\in\wU$ homogenous for the
bigradation $(\wU_{\nu}^{\delta})$
\begin{equation}\label{qscopcompmult}
\overline{\Delta}\bigl(\Psi(a)\Psi(b)\bigr)={(-1)}^{|a_2||b_1|} 
\overline{\Delta}\Psi(a)\overline{\Delta} \Psi(b)
\end{equation} 
with Sweedler notation
 $\Delta(a)=a_1\otimes a_2$, $\Delta(b)=b_1\otimes b_2$. 
Recall (\ref{qcompagrad}), (\ref{qcomcopbigrad})  which assert
in particular that $a,b,a_1,a_2,b_1,b_2$ are graded for the three
gradations  $(\wU_{\nu})$, $({}^{\mu}
\wU)$, $(\wU^{\delta})$.
By definition of $\overline{\Delta}$ and by (\ref{qcomPsi}) one has
$$\begin{array}{lcl}
\overline{\Delta}\bigl(\Psi(a)\Psi(b)\bigr)&=& 
{(-1)}^{(\nu(a),\mu(b)}\overline{\Delta}\Psi(ab)\\
&=& 
{(-1)}^{(\nu(a),\mu(b))+(\nu(a_1b_1),a_2b_2+\eta\nu(a_2b_2))}
\Psi(a_1b_1)\otimes \Psi(a_2b_2)\\
&=& {(-1)}^s \Psi(a_1)\Psi(b_1)\otimes \Psi(a_2)\Psi(b_2)
\end{array}$$
where $s\in \Z_2$ and 
$$s:=\bigl(\nu(a),\mu(b)\bigr)+
\bigl(\nu(a_1b_1),\nu(a_2b_2)+\eta\nu(a_2b_2)\bigr) +
\bigl(\nu(a_1),\mu(b_1)\bigr)+\bigl(\nu(a_2),\mu(b_2)\bigr).$$ 
According to (\ref{qcompagrad}) and
(\ref{qcomcopbigrad}), one has $\mu(b)=\delta(b)+\eta\nu(b)$, 
$\mu(b_1)=\delta(b)+\nu(b_2)+\eta\nu(b_1)$, $\mu(b_2)=\delta(b)+
\eta\nu (b_2)$, and hence
$$
\begin{array}{lcl}
s&=&\bigl(\nu(a_1a_2),\eta\nu (b_1b_2\bigr))+\bigl(\nu(a_1b_1),
\nu(a_2b_2)+ \eta\nu(a_2b_2)\bigr)+\bigl(\nu(a_1),\nu(b_2)\bigr)\\
&&\hskip 6.5truecm+\bigl(\nu(a_1),\eta\nu(b_1)\bigr)
+\bigl(\nu(a_2)+\eta\nu(b_2)\bigr).
\end{array}$$
Expending all scalar products in the above expression of $s$, we find
$$\begin{array}{lcl}s&=&
\bigl(\nu(a_1),\nu(a_2)+\eta\nu(a_2)\bigr)+
\bigl(\nu(b_1),\nu(b_2)+\eta\nu(b_2)\bigr)\\
&&\hskip 3truecm
+\bigl(\nu(a_2),\nu(b_1)\bigr)+\bigl(\nu(b_1),\eta\nu(a_2)\bigr)+
\bigl(\nu(a_2),\eta\nu(b_1)\bigr).
\end{array}$$
Using (\ref{qeqparity}), $s$ can be rewritten as
$$s= \bigl(\nu(a_1),\nu(a_2)+\eta\nu(a_2)\bigr)+\bigl(\nu(b_1),\nu(b_2)+
\eta\nu(b_2)\bigr)
+|a_2||b_1|$$
and therefore
$$\begin{array}{lcl}
\overline{\Delta}\bigl(\Psi(a)\Psi(b)\bigr)&=&
{(-1)}^{|a_2||b_1|}\Bigl({(-1)}^{(\nu(a_1),\nu(a_2)+\eta\nu(a_2))}
\Psi(a_1)\otimes \Psi(a_2)\Bigr)\\
&&\hskip 3.5truecm \times  \Bigl({(-1)}^{(\nu(b_1),\nu(b_2)+
\eta\nu(b_2))}  \Psi(b_1)\otimes \Psi(b_2)\Bigr)\\
&=&{(-1)}^{|a_2||b_1|}  \overline{\Delta}\Psi(a)\overline{\Delta} \Psi(b)
\end{array}$$
which proves (\ref{qscopcompmult}).

Let $m:\wU\otimes \wU\longrightarrow \wU$ be the multiplication map. 
We show next that for any element $a\in \wU$ homogenous for 
the bigradation $\wU_{\nu}^{\delta}$, 
\begin{equation}\label{qsantipodcompmult}
m(1\otimes \overline{S})\overline{\Delta}\Psi(a)=
m(\overline{S}\otimes 1)\overline{\Delta}\Psi(a)=\varepsilon\Psi(a)
\end{equation}
By definition of $\overline{\Delta}$, and $\overline{S}$
$$\begin{array}{lcl}
m(1\otimes \overline{S})\overline{\Delta}\Psi(a)&=&
{(-1)}^{(\nu(a_1),\nu(a_2)+\eta\nu(a_2))}\Psi(a_1)\otimes \overline{S}
\Psi(a_2)\\
&=&{(-1)}^{(\nu(a_1),\nu(a_2)+\eta\nu(a_2))+(\nu(a_2),\delta(a_2))}
\Psi(a_1)\otimes \Psi(Sa_2).
\end{array}$$
According to (\ref{qcomcopbigrad}), $S(a_2)\in \wU^{\nu(a_2)+
\delta(a_2)}_{\nu(a_1)}$ and $\delta(a_2)=\delta(a)$. Then formula
(\ref{qcomPsi}) gives
$\Psi(a_1)\Psi(Sa_2)={(-1)}^{(\nu(a_1),\delta(a)+\nu(a_2)+
\eta\nu(a_2))}\Psi(a_1Sa_2)$, and we obtain finally 
$$ m(1\otimes \overline{S})\overline{\Delta}\Psi(a)={(-1)}^{(\nu(a),
\delta(a))}\Psi(\varepsilon(a))=  {(-1)}^{(\nu(a),
\delta(a))}\varepsilon \Psi(a)=\varepsilon \Psi(a)$$ since
$\varepsilon \Psi(a)=0$ if $\nu(a)\not =0$. The other equality of 
(\ref{qsantipodcompmult}) can be established in the same way.

It remains to check that 
$$\begin{array}{c}
(\varepsilon \otimes 1)\overline{\Delta}\Psi(a)=(1\otimes \varepsilon)
\overline{\Delta}\Psi(a)=\Psi(a)\\
(\overline{\Delta}\otimes 1)\overline{\Delta}\Psi(a)=
(1\otimes\overline{\Delta})\overline{\Delta}\Psi(a).
\end{array}$$
These identities are straightforward from the definition of
$\overline{\Delta}$. 
\end{pf}

\section{Twisted adjoint actions and
locally finite parts}\label{qtwisad}
The object of this section is to compute the locally finite parts
of $\wU$ for certain twisted adjoint actions.

\subsection{A general construction}\label{qactadgenconstr}
 Let $X=X_{|_0}\oplus X_{|_1}$ be a Hopf superalgebra. 
We recall that the Hopf structure of $X$ provides an adjoint action
 defined by the formula 
$\ad a(x)={(-1)}^{|x||a_2|}a_1xS(a_2)$, using the Sweedler 
notation: $\Delta(a)=a_1\otimes a_2$. 
There is an elementary way to construct new actions by twisting the 
adjoint action by an algebra morphism (similar twisted  actions
have been  considered by Joseph in~\cite{jprim}).
One proceeds as follows. Let
$\psi:X\longrightarrow X$ be an algebra morphism. For
any $\Z_2$ homogenous element $a,x\in X$, set: 
$$\ad_{\psi}a(x)={(-1)}^{|a_2||x|}a_1xS(\psi(a_2)).$$
This formula defines an action  since for any homogenous elements $a,b,x$ 
$$\begin{array}{lcl}(\ad_{\psi}a)(\ad_{\psi}b)(x)&=&
{(-1)}^{|b_2||x|+|a_2|(|b_1|+|x|+|b_2|)}
a_1b_1xS(\psi(b_2))S(\psi(a_2))\\
&=&{(-1)}^{|x|(|a_2|+|b_2|)+|b_1||a_2|}a_1b_1xS(\psi(a_2)\psi(b_2))\\
&=&{(-1)}^{|x|(|a_2|+|b_2|)+|b_1||a_2|}a_1b_1xS(\psi(a_2b_2))\\
&=&{(-1)}^{|x||(ab)_2|}(ab)_1xS(\psi((ab)_2))\\
&=&\ad_{\psi}(ab)(x)
\end{array}$$

\subsection{}\label{qdeftwistadact} Let  $\mu\in
P_{\fg}(\pi)\slash 2P_{\fg}(\pi)$, and $\psi_{\mu}$ be the 
inner automorphism of $\wU$ defined by 
$\psi_{\mu}(a)=\xi_{\mu} a \xi_{\mu}$. 
In what follows we shall consider the twisted adjoint actions
obtained by applying the construction~\ref{qactadgenconstr}
to the cases of 
\begin{itemize}
\item the genuine Hopf algebra $\wU$ (for the Hopf structure given 
in~\ref{qHopf}) and the morphisms $\psi_{\mu}$.
\item the Hopf superalgebra $\wU$ (for the Hopf superstructure
given in~\ref{sHopfwU}) and the morphisms $\psi_{\mu}$.
\end{itemize}
In order to avoid any confusion, we shall write $\ad$  the adjoint
action of $\wU$, and $\bad$ the super adjoint action.
The twisted actions are denoted respectively by $\ad_{\mu}:=\ad_{\psi_{\mu}}$,
$\bad_{\mu}:=\bad_{\psi_{\mu}}$. 
Of course,  $\ad_{0}=\ad$ and $\bad_{0}=\bad$.
In the case $\mu=w_l$, we shall often prefer to write $\badp$ 
instead
of $\bad_{w_l}$. The twisted adjoint action $\badp$ 
is the quantum version of the ``non-standard'' 
adjoint action
 considered by Arnaudon, Bauer, Frappat in~\cite{abf}, 2. 
Recall~\ref{qdefsgen}.
By definition, $\badp a=\bad a$ for any generator 
$a\in \{k_{\lambda},\ \xi_{\lambda}  \lambda\in P_{\fg}(\pi);\
e_i, f_i, 1\leq i<l\}$ and 
\begin{equation}\label{eqadsdef}\begin{array}{c}
\bad(e_l)x=e_lx-{(-1)}^{|x|}k_lxk_l^{-1}e_l,\ \ 
\bad(f_l)x=f_lxk_l-{(-1)}^{|x|}xf_lk_l\\
\badp(e_l)x=e_lx+{(-1)}^{|x|}k_lxk_l^{-1}e_l,\ \ 
\badp(f_l)x=f_lxk_l+{(-1)}^{|x|}xf_lk_l
\end{array}\end{equation}

\subsubsection{}\label{qlocfindef} If $N$ is any $\ad_{\lambda}$-stable
(resp. $\bad_{\lambda}$-stable) subspace of $\wU$, we denote by 
$\F_{\lambda}(N)$ (resp. $\bF_{\lambda}(N)$) its locally finite part
for the action $\ad_{\lambda}$ (resp. $\bad_{\lambda}$). 
If $\mu=0$ we shall write respectively $\F(N)$, $\bF(N)$ instead of 
$\F_0(N)$, $\bF_0(N)$.   
Also, we shall often prefer to use
the notation $\bFp(N)$ instead of $\bF_{w_l}(N)$.

\subsubsection{}\label{qlemtwadact}
 Let $\lambda\in P_{\fg}(\pi)$. By definition of $\ad_{\lambda}$, \
for all $x\in \wU$, and for all weight
element $a\in \wU$, 
$(\ad_{\lambda} a)(\xi_{\lambda}x)=
a_1\xi_{\lambda}xS(\xi_{\lambda}a_2\xi_{\lambda})=
\xi_{\lambda} (\xi_{\lambda}a_1\xi_{\lambda}) x 
S(\xi_{\lambda}a_2\xi_{\lambda})={(-1)}^{(\lambda,\nu(a))}
\xi_{\lambda}(\ad a)x.$
The same holds replacing $\ad$ by $\bad$. 
It follows that 
$$\F_{\lambda}(\wU)=\xi_{\lambda} \F(\wU)\mbox{ and } 
\bF_{\lambda}(\wU)=\xi_{\lambda} \bF(\wU).$$

\subsubsection{} \begin{lem}{}\label{adactns0}
 Let $\lambda$ be in  $P_{\fg}(\pi)\slash
  2P_{\fg}(\pi)$. Then
\begin{enumerate} 
\item $\F_{\lambda}(\cU)=0$ if $\lambda\not =0$
\item $\F_{\lambda}(\wU)=\xi_{\lambda}\F(\cU)$
\end{enumerate}
\end{lem}
\begin{pf} Assume that $\lambda\not=0$, and let $V\subset \cU$ be a simple
$\ad_{\lambda}\cU$-module. Take $a$ an element of lowest weight of $V$.
Since $\lambda\not=0$, there exists
$\alpha_i\in \pi$ such that $(\lambda,\alpha_i)=1+2\Z$.
Hence $0=\ad_{\mu} F_ia=(F_ia+aF_i)K_{\alpha_i}$.
Proposition 1.7,~\cite{dk}, forces $a=0$. This establishes the
assertion (i). 

Let $a\in \F(\wU)$. Write $x=\sum_{\mu\in
P_{\fg}(\pi)\slash 2P_{\fg}(\pi)} \xi_{\mu}a_{\mu}$, $x_{\mu}\in 
\cU$.
According to~\ref{qlemtwadact}
$\ad (\wU)x=\ad(\cU)x=
\bigoplus_{\mu\in
P_{\fg}(\pi)\slash 2P_{\fg}(\pi)} \xi_{\mu}\ad_{\mu}(\cU)(x_{\mu}).$
Thus  $\ad_{\mu}(\cU)(x_{\mu})\in \F_{\mu}(\cU)$ and hence
$x_{\mu}=0$ if $\mu\not=0$ by (i). This proves $\F(\wU)=\F(\cU)$ and
(ii) follows from~\ref{qlemtwadact}. 
\end{pf}

\subsection{} Retain the definitions of~\ref{qdefgraddelta}.
If follows from (\ref{qcomcopbigrad}) and (\ref{qeqdefsHopf}) that
the gradation $(\wU^{\delta})$ possesses a very striking property:
it is invariant by the actions $\ad_{\lambda}, \bad_{\lambda}$.
Recall the definition of $\Psi$ (see~\ref{defzhtransfor}). 

\begin{lem}{}\label{qlemtransfadbad}  
Fix $\lambda,\ \delta\in {P_{\fg}(\pi)/2P_{\fg}(\pi)}$.
Let $a,x\in \wU$ be 
homogenous elements for the bigradation $(\wU^{\delta}_{\nu})$. 
Then
\begin{equation}\label{formadadb}
\Psi\bigl(\ad_{\lambda}a(x)\bigr)=\pm
\bad_{\lambda+\delta}\Psi(a)\bigl(\Psi(x)\bigr).
\end{equation}
\end{lem}
\begin{pf} Let $a,x\in \wU$ be as in the Lemma. 
By definition, $\ad_{\lambda}a x={(-1)}^{(\nu(a_2),\lambda)}
a_1 x S(a_2)$,  where $\Delta(a)=a_1\otimes a_2$ in Sweedler notation.
Recall (\ref{qcompagrad}), (\ref{qcomcopbigrad})  which assert
in particular that $a,x,a_1,a_2$ are graded for the
gradations  $(\wU_{\nu})$, $({}^{\mu}
\wU)$, $(\wU^{\delta})$.
Using (\ref{qcomPsi}) one obtains
\begin{align}
\Psi(\ad_{\lambda}a x)&={(-1)}^{(\nu(a_2),\lambda)+(\nu(a_1),\mu(x))+
(\nu(x),\mu (Sa_2))+(\nu(a_1),\mu (Sa_2))}
\Psi(a_1)\Psi(x)\Psi(Sa_2)\nonumber\\
&= {(-1)}^s \Psi(a_1)\Psi(x)\overline{S}\Psi(a_2)\label{qeqadcom1}
\end{align}
where $s\in \Z_2$, and 
$$s=\underbrace{\bigl(\nu(a_2),\lambda\bigr)+\bigl(\nu(a_1),\mu(x)\bigr)+
\bigl(\nu(x),\mu (Sa_2)\bigr)}_{s_2}+
\underbrace{\bigl(\nu(a_1),\mu (Sa_2)\bigr)+\bigl(\nu(a_2),\delta(a_2)
\bigr)}_{s_2}.$$
According to (\ref{qcomcopbigrad}) $\mu(x)=\delta(x)+\eta\nu(x)$, 
$\delta(a_2)=\delta(a)$, $\mu(S(a_2))=\delta(a)+\nu(a_2)+\eta\nu(a_2)$,
so
$$s_1=\bigl(\nu(a),\delta(a)\bigr)+\bigl(\nu(a_1),\nu(a_2)+\eta\nu(a_2)
\bigr)$$
and 
$$\begin{array}{lcl}
s_2&=& \bigl(\nu(a_2),\lambda\bigr)+\bigl(\nu(a_2)+\nu(a),\delta(x)+
\eta\nu(x)\bigr)
+\bigl(\nu(x),\nu(a_2)+\eta\nu(a_2)+\delta(a)\bigr)\\
&=& \bigl(\nu(a),\mu(x)\bigr)+\bigl(\nu(x),\delta(a)\bigr)+\bigl(\nu(a_2),
\lambda+\delta(x)\bigr)\\
&& \hfill+ \bigl(\nu(a_2),\eta\nu(x)\bigr)+
\bigl(\nu(x),\eta\nu(a_2)\bigr)+\bigl(\nu(x),\nu(a_2)\bigr)\\
&=& \bigl(\nu(a),\mu(x)\bigr)+\bigl(\nu(x),\delta(a)\bigr)+\bigl(\nu(a_2),
\lambda+\delta(x)\bigr)+
|a_2||x|\ \ \mbox{ by (\ref{qeqparity})}.
\end{array}$$
Consequently $$s=s_1+s_2=t + |a_2||x| + 
(\nu(a_2),\lambda+\delta(x))+\bigl(\nu(a_1),\nu(a_2)+\eta\nu(a_2)\bigr)$$
where $t=(\nu(a),\delta(a))+ (\nu(a),\mu(x))+(\nu(x),\delta(a))$ 
depends only on the ``degrees'' (for the different gradations)
 of $a$ and $x$. Substituting the above expression of $s$ in 
(\ref{qeqadcom1}), and using  definitions of 
$\overline{\Delta}$, $\overline{S}$ we derive  that 
\begin{equation}\label{qeqfinecompad}
\Psi(\ad_{\lambda}a x)={(-1)}^t \bad_{\lambda+\delta(x)}\Psi(a)\Psi(x)
\end{equation} 
as required.  
\end{pf}

\subsection{}\label{qalgstrucrec}
 We recall the results of Joseph and Letzter (see~\cite{jl2}):
$$\F(\cU)=\displaystyle\bigoplus_{\lambda\in P^+_{\fk}(\pi)} (\ad\cU) 
K_{-2\lambda}$$
and each $(\ad\cU) K_{-2\lambda}$ contains a unique (up to a non-zero 
scalar) central element denoted by $z_{2\lambda}$. The centre $\cZ(\cU)$ of 
$\cU$ is the polynomial algebra 
\begin{equation}\label{tonycentre}
\cZ(\cU)=\C[z_{2w_1},\ldots, z_{2w_{l-1}}, z_{w_l}].
\end{equation} 
We shall need the following submodules of $\F(\cU)$:
$$N_0:=\displaystyle\bigoplus_{\lambda\in P^+_{\fg}(\pi)}
(\ad\cU)K_{-2\lambda},\ \ \ N_1:=
\displaystyle\bigoplus_{\lambda\in P^+_{\fk}(\pi)
\backslash P^+_{\fg}(\pi)} (\ad\cU) K_{-2\lambda}.$$
If $\lambda\in P^+_{\fg}(\pi)$, then  
$\delta(K_{-2\lambda})=2\lambda=0$ and so $N_0\subset \cU^0$.
If ${\lambda\in P^+_{\fk}(\pi)
\backslash P^+_{\fg}(\pi)}$, i.e. 
$\lambda=\lambda'+{w_l\over 2}$, $\lambda'\in   P^+_{\fg}(\pi)$, then 
$\delta(K_{-2\lambda})=w_l$ and so  $N_1\subset \cU^{w_l}$.
Consequently 
\begin{equation}\label{qN0N1}
\begin{array} {l} \F(\cU)\cap \cU^0=N_0\\ 
\F(\cU)\cap \cU^{w_l}=N_1.
\end{array}
\end{equation}
It follows from~\Lem{qlemtransfadbad}  that
$\bF(\wU)\cap \wU^{\mu}=\Psi\bigl(\F_{\mu}(\wU)\cap \wU^{\mu}\bigr)$.
By~\Lem{adactns0} we know that $\F_{\mu}(\wU)=\xi_{\mu}\F(\cU)$. 
Hence, using (\ref{qN0N1}) and recalling~\ref{qlemtwadact}, we obtain
$\forall \lambda\in P_{\fg}(\pi)/2 P_{\fg}(\pi)$,
\begin{equation}\label{qeqbFfctF2}
\begin{array}{lcl}
\bF_{\lambda}(\wU)&=&\bigoplus_{\mu} \xi_{\lambda+\mu}\Psi\bigl(\F(\cU)
\cap \cU^{\mu}\bigr)\\
&=&\xi_{\lambda}\bigl( \Psi(N_0)\oplus \xi \Psi(N_1)\bigr).
\end{array}
\end{equation}

\section{Algebraic structures of $\cU_{q}(\fg)$}\label{secalgstrucbU}

Recall the definition of $\Psi$ (see~\ref{defzhtransfor}).
We define $\bU:=\Psi(\cU)$. By definition of $\Psi$,
$$\wU\simeq \bU\rtimes k\Gamma.$$
With the notations of~\ref{qdefsgen}, $\bU$ is
 the subalgebra of $\wU$ generated by  the $e_i$, $f_i$, 
$k_{\lambda}$. The algebra $\bU$ is graded for all the different
gradations we defined on $\wU$. 
By definition of $\overline{\Delta}, \overline{S}$,
the subalgebra 
$\bU$ is a Hopf subalgebra of $(\wU,\overline{\Delta}, \overline{S},
\varepsilon)$.

We shall show in~\ref{qalgstUs} that $\bU\simeq \cU_{-q}(\fg)$.

\subsection{Representations of $\bU$}
We shall now give the classification of the
finite dimensional $\bU$-modules.

 \subsubsection{Generalities}\label{fdUsmod.1}
Let us prove that every simple $\bU$-module is
the restriction of a simple $\wU$-module and
that every finite dimensional $\bU$-modules
is completely reducible.

Let $V$ be a simple finite dimensional $\bU$-module.  Assume for the moment
that we work over the algebraic closure $\overline{k}$ of $k$ (we
extend the scalars of all our objects).  
Take any non-zero weight vector of $V$ (that is a
common eigenvector for the $k_{\lambda}$).
The simplicity forces $V=\bU v$. Choose any character 
$\theta\in \Gamma$. Then the the following 
 formula defines an action of $\wU$ on
$V$ : $\xi_{\mu}.{a_{\nu}v}=\theta(\mu){(-1)}^{(\mu,\nu)}a_{\nu}v$,
 for any $\nu\in P_{\fg}(\pi)$ and any
$a_{\nu}\in \bU$ of weight $\nu$. Indeed, the vector $v$ being
a weight vector, the annihilator $\Ann_{\bU} v$ is the sum
of its weight subspaces, and hence the previous formula makes sense.
As a $\wU$-module,
$V$ is necessarily simple, and thus is a 
$V(\phi q^{\lambda},\theta)$ by  lemma~\ref{propcomplredwU} (which
obviously also holds  over $\overline{k}$). This shows in particular
that all the eighenvalues of the $k_{\lambda}$ actually lie in our 
ground field $k$. Therefore, we could have chosen $v$ to be 
in the $k$-vector space $V$, and so all that precedes actually 
holds over $k$. We have proved that $V$ is the restriction of a
simple $\wU$-module.

Remark that we have just showed that the group 
$\{k_{\mu},\ \mu\in P_{\fg}(\pi)\}$ acts diagonally on a simple 
$\bU$-module.
Hence the restriction of a simple finite dimensional $\wU$-module  to
$\bU$
is also simple.

Consider now  $M$, a finite dimensional 
$\bU$-module. The induced
module $\Ind_{\bU}^{\wU}M$ is a finite dimensional $\wU$-module, and 
therefore completely reducible by lemma~\ref{propcomplredwU}. Thus,
$\Ind_{\bU}^{\wU}M$ is also completely reducible as a $\bU$-module (see
the previous remark). 
 But as a $\bU$-module,
$M$ lies in $\Ind_{\bU}^{\wU}M$. Hence 
$M$ is completely reducible.

\subsubsection{} It follows from what preceedes,
that for any fixed $\theta\in\widehat{\Gamma}$, 
the set $\{V(q^{\lambda}\phi,\theta),\ 
(\lambda,\phi)\in P_{\fk}^+(\pi)\times\widehat{\Gamma}\}$ is a complete set 
 of non-isomorphic 
finite dimensional simple modules for both $\cU$ and $\bU$.

\subsubsection{}\begin{rem}{}  Recall that we shall prove in~\ref{qalgstUs} 
that $\bU\simeq \cU_{-q}(\fg)$.
The classification of the finite dimensional modules over the 
``quantum'' enveloping algebra of $\fg$ has been 
obtained by R.~B.~Zhang (see~\cite{zh}) 
in the context of formal deformations
 and by  Zou (see~\cite{zo}) for the Drinfeld-Jimbo quantization
$\cU_q(\fg)$,
 through the standard approach.
\end{rem} 

\subsubsection{Crystals}
We admit for a moment 
that $\bU\simeq\cU_{-q}(\fg)$.
Fix $\theta\in \widehat{\Gamma}$ and 
let $V(\lambda)$, $\lambda\in P_{\fk}^+(\pi)$, be the simple finite
dimensional $\wU$-module $V(q^{\lambda},\theta)$.
{\it A priori}, one can associate to $V(\lambda)$ two crystals. One is given
by the work of Kashiwara (see~\cite{kash}), considering $V(\lambda)$ as
a $\cU$-module. We denoted it by $B(\lambda)$. The  other one, 
${\cal B}(\lambda)$, follows from the work of 
 Musson and Zou (see~\cite{mz}),
viewing this time $V(\lambda)$ as a $\bU$-module. 
Both sets $\{B(\lambda),\ \lambda\in P_{\fk}^+(\pi)\}$,
$\{{\cal B}(\lambda),\ \lambda\in P_{\fk}^+(\pi)\}$ are closed family
of highest weight normal crystals, in the sense of~\cite{jq} 6.4.21.
Hence  there are equal (up to isomorphisms) by 
 proposition 6.4.21,~\cite{jq}.
Another way to see that $B(\lambda)\simeq {\cal B}(\lambda)$ is to 
remark (keeping the notations of~\cite{kash} and~\cite{mz}) that
${\cal L}(\lambda)=L(\lambda)$ and that the crystalline operators 
of Musson and Zou act on a given weight subspace of $L(\lambda)$
as the crystalline operators of Kashiwara up to signs 
(depending on 
the weight of the subspace 
and on the ``color'' of the operators).

\subsection{} Recall the definition of $\cU_q(\fg)$ given in~\cite{mz}.
Let $\bZ(\bU)$ (resp. $\bZp(\bU)$) be the supercentre 
(resp. the anticentre) 
of $\bU$, that is
the subspace of invariants elements of $\bU$ with respect to $\bad$
(resp $\badp$). One has 
$\bZp(\bU):=\{a\in \bU,\ ax={(-1)}^{|x|}xa\ \forall \ \Z_2
\mbox{-homogenous } x\in \bU\}$.  We also introduce  the 
algebra $\cA(\bU):=\bZ(\bU)+\bZp(\bU)$.

We deduce from section~\ref{qtwisad} the (compare (i) with 3.3~\cite{zh})

\begin{thm}{}\label{qalgstUs}  \begin{enumerate}
\item The subalgebra $\bU$ is isomorphic to $\cU_{-q}(\fg)$.
\item One has $\Psi\bigl(\F(\cU)\bigr)=\bF(\bU)\oplus \bFp(\bU)$ with
$$\begin{array}{l}
\bF(\bU)=\Psi(N_0)= \displaystyle\bigoplus_{\lambda\in P^+_{\fg}(\pi)}
\bad_{}\bU k_{-2\lambda}\end{array}$$$$
\begin{array}{l}
 \bFp(\bU)=\Psi(N_1)=\displaystyle\bigoplus_{\lambda\in P^+_{\fk}(\pi)
\backslash P^+_{\fg}(\pi)}\bad_{}\bU
k_{-2\lambda}\end{array}$$
\item Recall that $\xi:=\xi_{\beta_l}$
and (\ref{tonycentre}).  One has $\cA(\bU)=\bZ(\bU)\oplus
\bZp(\bU)=\Psi(\cZ(\cU))$ and
$$\begin{array}{l}
\bZ(\bU)=\C[z_{2w_1},\ldots,z_{2w_{l-1}},z_{w_l}^2]\\
\bZp(\bU)=(\xi z_{w_l})\bZ(\bU)
\end{array}$$
\end{enumerate}
\end{thm}
\begin{pf} We start by proving (i).  The relations (\ref{defwueq3}),
(\ref{defwueq4}), (\ref{reluhsere}), (\ref{reluhserf}) can be 
respectively rewritten as $\ad K_{\mu}E_i=q^{(\mu,\alpha_i)}E_i$,
$\ad K_{\mu}F_i=q^{-(\mu,\alpha_i)}F_i$, $\ad F_iE_j =\delta_{ij}
(1-K_{\alpha_i}^2)/ (q-q^{-1})$, $\ad E_i^{1-\langle \alpha_j,
\alpha_i\rangle} E_j=0$,  $\ad F_i^{1-\langle \alpha_j,
\alpha_i\rangle} F_j=0$. Take the image of these relations by $\Psi$.
According to Lemma~\ref{qlemtransfadbad}  
 (and more precisely to formula (\ref{qeqfinecompad})) we
obtain $\bad k_{\mu}e_i={(-q)}^{(\mu,\alpha_i)}e_i$, 
$\bad k_{\mu}f_i={(-q)}^{-(\mu,\alpha_i)}f_i$,
$\bad f_ie_j =-\delta_{ij}{(-1)}^{\delta_{il}}
(1-k_{\alpha_i}^2)/ (q-q^{-1})$, $\ad e_i^{1-\langle \alpha_j,
\alpha_i\rangle} e_j=0$, $\ad f_i^{1-\langle \alpha_j,
\alpha_i\rangle} f_j=0$. It is easy to see that 
these relations are exactly the relations defining 
$\cU_{-q}(\fg)$.

The assertion (ii) results from (\ref{qeqbFfctF2}). And (ii) implies
that
 $\Psi(\cZ(\cU)\cap N_0)=\bZ(\bU)$ and $\Psi(\cZ(\cU)\cap N_1)=\bZp(\bU)$.
On the other hand, elements of the centre $\cZ(\cU)$ 
are of weight zero. Therefore
combining (\ref{qN0N1}) with (\ref{qcompagrad}),
one has $\Psi(\cZ(\cU)\cap N_0)=\cZ(\cU)\cap N_0$ and
$\Psi(\cZ(\cU)\cap N_1)=\xi (\cZ(\cU)\cap N_1)$,
which ends the proof.
\end{pf}

\subsubsection{}\begin{rem}{}
 The element $\xi z_{w_l}$ is a quantization of the
element $T$ introduced in 4.4.1~\cite{gl2}. 
See also formula (\ref{qformhcT}) and 
remark~\ref{qpropqT}. Notice also that this element coincides with
the sCasimir element constructed in~\cite{ab} for the algebra
$\cU_q(\osp(1,2))$
\end{rem}

\subsubsection{}\begin{rem}{} Let $\bFp(\cU(\fg))$ be the 
locally finite part of $\cU(\fg)$ for the ``non-standard'' 
action of Arnaudon, Bauer, Frappat 
(see~\cite{abf}, 2). Then 
$\bFp(\cU(\fg)))=\bF(\cU(\fg))$ since $\fg_1$
is finite dimensional. The situation in the quantum case is therefore
radically different on this point. 
\end{rem}

\subsection{The separation theorem for $\bF(\bU)\oplus \bFp(\bU)$} 
\label{qsecsepthbU}
In~\cite{jl3}, Joseph and Letzter established a separation theorem for
the algebra $\F(\cU)$. They proved the existence of 
$\ad$-submodules $\cH(\cU)(\lambda)$ of $(\ad\cU) K_{-2\lambda}$
such that if $\cH(\cU):=\oplus_{\lambda\in
  P^+_{\fk}(\pi)}\cH(\cU)(\lambda)$,
then the multiplication 
$\cH(\cU)\otimes \cZ(\cU)\longrightarrow \cU$ is an isomorphism of
$\ad \cU$-modules. Introduce 
$$\bH(\bU):= \Psi\bigl(N_0\cap \cH(\cU)\bigr) 
\mbox{ and }\bHp(\bU):=\Psi\bigl(N_1\cap \cH(\cU)\bigr) $$
Let $h_i\in \cH(\cU)$ be weight elements, and 
$z_i\in \cZ(\cU)\cap N_{j_i}$ $j_i\in \{0,1\}$.
 It follows from (\ref{qN0N1}) and (\ref{qcompagrad})  that 
$\sum \Psi(h_i)\Psi(z_i)=\Psi(\sum \pm h_iz_i)$. This is enough to 
deduce the separation theorem for $\bF(\bU)\oplus \bFp(\bU)$:
\subsubsection{}\label{qsepthUs} 
\begin{prop}{}The multiplication $$(\bH(\bU)\oplus \bHp(\bU))\otimes
\cA(\bU)\longrightarrow \bF(\bU)\oplus \bFp(\bU)$$
is an isomorphism.
\end{prop}

\section{The annihilation theorem}\label{qsecdufloth}

The goal of this section (theorem~\ref {qthmduflo}) is to establish that 
the annihilator of any $\bU$-Verma module in $\bF(\bU)\oplus\bFp(\bU)$
is generated by its intersection with $\cA(\bU)$.

\subsection{}\label{qgenopVmod}  By definition (see~\ref{qrepgen}), a
$\wU$-Verma module is a $\cU$-Verma module and 
a character of 
$\Gamma$ which describes the action of $\Gamma$
 on the highest weight vector.    
Of course, the same holds replacing $\cU$ by  $\bU$. Hence,
throughout this subsection, 
we shall not  make any distinctions
between the $\cU$, $\bU$ and $\wU$-Verma modules.

We fix now once for all a 
$\wU$-Verma module $M:=M(\Lambda,\theta)$ and we define for all 
$\mu\in P_{\fg}(\pi)$, $\Lambda(k_{\mu}):=
\Lambda(K_{\mu})\theta(\xi_{\mu})$.

\subsection{}\label{qsecdeflocaendVm}
The Verma module $M$ being $\Z_2$-graded (see~\ref{qrepgengradM}),
 $\End(M)$ inherits  
the natural gradation: 
$${\End(M)}_{|_{{j}}}=
\{f\in \End(M), \forall i\in\Z_2\  f(M_{|_{ {i}}})\subset 
M_{|_{{i}+{j}}}\}$$

Consider the adjoint action of $\wU$ on  $\End(M)$ defined by 
$(\ad a f)(x)=a_1
(f({S}(a_2)x))$ for all $a\in \wU,\ f\in \End(M)$
and all $x\in M$. 
Let $\F(M,M)$ be the locally finite part of $\End(M)$
for the adjoint action $\ad$. The subspace $\F(M,M)$ is 
$\Z_2$-graded for the above gradation. 
The restriction of 
$\wU\longrightarrow \End(M)$ induces a   morphism
of $\ad\wU$-modules: 
$\F(\wU)\longrightarrow \F(M,M)$. Its image coincides
with the image of $\F(\cU)\longrightarrow \F(M,M)$. 

We recall (see Lemma 8.3. in~\cite{jq}) that $\F(M,M)$ is a domain.

\subsection{}
\begin{lem}{}\label{qinteropVm} Let $f\in \F(M,M)$ and $i\in \Z_2$
be be such that $f(M_{|_i})=0$. Then $f=0$.
\end{lem}
\begin{pf} Let $f$ be as in the lemma. 
By definition of the $\Z_2$-gradation on $\F(M,M)$ we may
assume that $f$ is $\Z_2$-homogenous, and hence that $f^2$ is even.
 Take any non-zero $p\in  {\F(M,M)}_{|_1}$ (obviously such $p$
exists; for instance $\ad E_l K_{-w_l}=(1-q^{-1})E_lK_{-w_l}\in \F(\cU)$ 
has a non-trivial
image in ${\F(M,M)}_{|_1}$).
Then $f^2pf^2=0$  which implies $f=0$ since $\F(M,M)$ is a domain.
\end{pf}

\subsection{} Recall that 
$\Lambda(k_{\mu}):=\Lambda(K_{\mu})\theta(\xi_{\mu})$.

\begin{lem}{}\label{qlemqduflo} For any $\bU$-Verma module $M$,
the following equivalence holds
$$\Ann_{\cA(\bU)}M=\cA(\bU)\Ann_{\cZ(\bU)}M\Longleftrightarrow
 \forall \ 1\leq i\leq l,\ \Lambda(k_{\beta_i})\not=
\pm iq^{-(\rho,\beta_i)}.$$
\end{lem} 
\begin{pf} By proposition~\ref{qalgstUs}, one has 
$\cA(\bU)=(k\oplus k(\xi z_{w_l}))\otimes \cZ(\bU)$. The
centre $\cZ(\bU)$
acts by scalars on $M$.  Hence  the equality
$\Ann_{\cA(\bU)}M=\cA(\bU)\Ann_{\cZ(\bU)}M$ is equivalent to 
$\Ann_{k\oplus k(\xi z_{w_l})}M=0$. The element 
$\xi z_{w_l}$ acts on the $\Z_2$-graded components of $M$ by the two  
opposite scalars $\pm\Lambda(\Upsilon(z_{w_l}))$, 
 $\Lambda$ being linearly extended 
to $\cU^o$. 
Thus, $\Ann_{k\oplus k(\xi z_{w_l})}M=0$ is equivalent to 
$\Lambda(\Upsilon(z_{w_l}))\not=0$. Retain notation of~\ref{defhcq}.
By~\cite{jq} 7.1.19,
\begin{equation}\label{qformhcz}
\Upsilon (z_{w_l})=\displaystyle \sum_{\mu\in P_{\fk}(\pi)}
\dim {V(q^{{w_l\over 2}})}_{\mu} q^{-2(\rho, \mu)}K_{-2\mu}
\end{equation} where 
$V(q^{{w_l\over 2}})$ is simple $\cU$-module of highest weight
${w_l\over 2}$.  Since ${w_l\over 2}$ is a minuscule weight, 
$\dim {V(q^{{w_l\over 2}})}_{\mu}=1$ if 
$\mu\in W ({w_l\over 2})=
\{ {1\over 2}\sum_{i=1\ldots l}
\varepsilon_i\beta_i,\ \varepsilon_i=\pm 1\}$ and
$0$ otherwise. So (\ref{qformhcz}) can be rewritten
\begin{equation}\label{qformhcT}
\Upsilon (z_{w_l})=\displaystyle \prod_{i=1,\ldots l}
(q^{-(\rho,\beta_i)}K_{-\beta_i}+q^{(\rho,\beta_i)}K_{\beta_i}).
\end{equation}
The assertion follows.
\end{pf}

\subsubsection{} \begin{rem}{}\label{qpropqT} If 
$\widehat{T}_d:=\{\Lambda\in\widehat{T},\ \exists 1\leq i\leq t
\mbox { such that } \Lambda(K_{\beta_i})=\pm iq^{-(\rho,\beta_i)}\}$
then the formula (\ref{qformhcT})  implies that
$$\xi z_{w_l}\in \bigcap_{\Lambda\in \widehat{T}_d,\ \theta\in 
\widehat{\Gamma}}
\Ann_{\bU}M (\Lambda,\theta). $$
This is the quantum version of a property satisfied by the element
$T$ constructed in~\cite{gl2} (see~\cite{gl2} 4.4.1 and 6.1.3).
\end{rem}

\subsection{} In ~\cite{jl}, Joseph and Letzter prove the
annihilation
theorem of Duflo for $\F(\cU)$. 
We  deduce from this result the (compare with Theorem 7.1~\cite{gl},
and Theorem 6.2~\cite{gl2}) 

\begin{thm}{}\label{qthmduflo} Let $M$ be a $\bU$-Verma module. Then
\begin{enumerate}
\item 
For any $i\in \Z_2$, 
$\Ann_{\bF(\bU)\oplus \bFp(\bU)}M_{|_i}=(\bF(\bU)\oplus \bFp(\bU))
\Ann_{\cA(\bU)}M_{|_i}.$
\item $\Ann_{\bF(\bU)\oplus \bFp(\bU)}M=(\bF(\bU)\oplus \bFp(\bU))
\Ann_{\cA(\bU)}M.$
\item $\Ann_{\bF(\bU)\oplus 
\bFp(\bU)}M=(\bF(\bU)\oplus \bFp(\bU))
\Ann_{\cZ(\bU)}M \Longleftrightarrow 
\forall \ 1\leq i\leq l,\ \Lambda(k_{\beta_i})\not=\pm iq^{-(\rho,\beta_i)}$
\end{enumerate}
\end{thm}
\begin{pf} 
We start by (i). Recall~\Thm{qalgstUs} (iii).  
As the centre $\cZ(\bU)$ acts by scalars on 
$M$, the algebra $\cA(\bU)$ acts by scalars on the homogenous
components $M_{|_{{i}}}$. It follows from 
Proposition~\ref{qsepthUs}  that (i) 
is equivalent to the statement 
$\forall i\in\Z_2$, $\Ann_{\bH(\bU)\oplus \bHp(\bU)}M_{|_i}=0$.
Let
$\Psi(h),\Psi(h')\in \bH(\bU),\bHp(\bU)$ and $i\in \Z_2$
 be such that
$\Psi(h)+\Psi(h')\in\Ann_{\bH(\bU)\oplus \bHp(\bU)}M_{|_i}$.
Since $M_{|_i}$ is  $T$-invariant,
 we can assume that $\Psi(h),\Psi(h')$ 
(and hence $h,h'$) are elements of the same weight $\nu$. 
Combining (\ref{qN0N1}) and (\ref{qcompagrad}) 
one has  $\Psi(h)=\xi_{\eta(\nu)} h$ and 
$\Psi(h')=\xi_{\eta(\nu)+w_l} h'$.
Hence $h +\xi h'\in \Ann_{\cH(\cU)} M_{|_i}$. The element $\xi$
acts by $\pm id$ on $M_{|_i}$, so we can assume that
$h + h'\in \Ann_{\cH(\cU)} M_{|_i}$.
From Lemma~\ref{qinteropVm} we derive that 
$h+h'\in \Ann M$.
Therefore $h'=-h$  using 4.2,~\cite{jl}. But $h\in N_0$,
$h'\in N_1$ and $N_0\cap N_1=0$, which forces finally 
$h=h'=0$. This finishes the proof of (i).

For (ii), one has the equalities
$$\begin{array}{lcl} \Ann_{\bF(\bU)\oplus 
\bFp(\bU)}M&=& \bigcap_i \Ann_{\bF(\bU)\oplus 
\bFp(\bU)}M_{|_i}\\
&=& \bigcap_i (\bF(\bU)\oplus \bFp(\bU))
\Ann_{\cA(\bU)}M_{|_i}\ \mbox{ by (i)}\\
&=& \bigcap_i (\bH(\bU)\oplus \bHp(\bU))\otimes 
\Ann_{\cA(\bU)}M_{|_i}\ \mbox { by~\Prop{qsepthUs}} \\
&=&(\bH(\bU)\oplus \bHp(\bU))\otimes  \bigcap_i \Ann_{\cA(\bU)}M_{|_i}\\
&=&(\bH(\bU)\oplus \bHp(\bU))\otimes \Ann_{\cA(\bU)} M\\
&=& (\bF(\bU)\oplus \bFp(\bU))
\Ann_{\cA(\bU)}M
\end{array}$$
And (iii) is a consequence of (ii),~\Prop{qsepthUs} and of 
Lemma~\ref{qlemqduflo}.
\end{pf}

\subsubsection{}\begin{rem}{} We believe that (i) should also hold in 
the classical case.
\end{rem}

\section{$\cA(\bU)$ is the commutant of ${\bU}_{|_0}$ }
\label{qseccom}

In this section we shall prove
that $\cA(\bU)$ is the commutant of the even part of $\bU$, that is 
$\cA(\bU)=\cC(\bU_{|_0})$. Let $\cA(\wU)$ be the subalgebra
 $\cA(\wU):=\cZ(\cU)\oplus \xi \cZ(\cU)$. Since
 $\cC({\wU}_{|_0})\cap \bU=\cC(\bU_{|_0})$ and 
$\cA(\wU)\cap \bU=\cA(\bU)$ (recall~\Thm{qalgstUs} (iii))
it is enough to prove the
equality $\cA(\wU)=\cC({\wU}_{|_0})$. For this, we shall proceed 
by quantizing the mechanics of 4~\cite{gl2}.

\subsection{}\label{interannusimq}
\begin{lem}{} For every subset $\Omega$ of $P_{\fk}(\pi)$ dense for
the Zariski topology, one has
 $$\bigcap_{\lambda\in\Omega\atop \theta\in\widehat{\Gamma}}
\Ann_{\wU} V(q^{\lambda},\theta)=0.$$
\end{lem}
\begin{pf} Fix $\lambda\in\Omega$. One has 
\begin{equation}\label{eqannusm}
\bigcap_{\theta\in\widehat{\Gamma}}\Ann_{\wU} V(q^{\lambda},\theta)
=(k\Gamma)\Ann_{\cU} V(q^{\lambda}).
\end{equation}
where $ V(q^{\lambda})$ stands for the $\cU$-simple module of 
highest weight $q^{\lambda}$.

Indeed, let $e_{\chi}\in k\Gamma,\ \chi\in\Gamma$, be the projector
corresponding to $\chi$, that is the projector
such that $ge_{\chi}=\chi(g)e_{\chi}$,
$\forall g\in \Gamma$. Let $x\in\wU$ and write 
$x=\sum_{\chi\in \widehat{\Gamma}}x_{\chi}e_{\chi}$, $x_{\chi}\in \cU$.
As a $\cU$-module, $V(q^{\lambda},\theta)$ canonically identifies with the 
$\cU$-module $V(q^{\lambda})$  (see~\ref{qrepgen}). Under this 
identification,
 $x$ acts on 
the subspace of $T$-weight $q^{\lambda-\nu}$  of $V(q^{\lambda},\theta)$ 
as $x_{{(-1)}^{\nu}\theta}$ on the subspace of the same weight of
 $V(q^{\lambda})$. Hence, $x\in   
\bigcap_{\theta\in\widehat{\Gamma}}\Ann_{\wU}
  V(q^{\lambda},\theta)$ implies that $x_{{(-1)}^{\nu}\theta}$ vanishes
on $V(q^{\lambda})_{q^{\lambda-\nu}}$ for all $\nu$ and $\theta$.
This gives (\ref{eqannusm}).
\end{pf}

\subsection{}  As in 4.1~\cite{gl2}, the previous 
lemma implies 

\begin{lem}{}\label{qcommUactscal} The algebra  $\cC({\wU}_{|_0})$ 
coincides with  the subalgebra of 
elements of $\wU$ acting by scalars on the homogenous components of simple
highest weight modules.
\end{lem}
\subsection{}\label{qpfcentnotchar} Retain the notation of~\ref{defhcq}.
Take $x\in \wU^o$ and write $x:=\sum a_{\mu,\mu'}\xi_{\mu}K_{\mu'}$,
$a_{\mu,\mu'}\in k$.
For any $(\lambda,\theta)\in P_{\fk}(\pi)\times
\widehat{\Gamma}$, we set $x(\lambda,\theta):=\sum a_{\mu,\mu'}
\theta(\mu)q^{(\lambda,\mu')}$.

With these conventions, $a\in\cC(\wU_{|_0})$ acts by the scalar 
$\Upsilon(a)(\lambda,\theta)$ on the homogenous component of 
$V(q^{\lambda},\theta)$ containing the highest weight vector. 

\subsection{}
\begin{lem}{} \label{lemhccentwu1}
 The restriction of  $\Upsilon$ to $\cC({\wU}_{|_0})$ is
  injective.
\end{lem}
\begin{pf}  For all $(\lambda,\theta)\in P^+_{\fk}(\pi)\times
\widehat{\Gamma}$ denote by  $v_{\lambda}$ 
the highest weight vector of  $V(q^{\lambda},\theta)$. Let
$a$ be in $\cA(\wU)$. We recall that $a$ acts on
$\wU_{|_0}v_{\lambda}$
by the scalar $\Upsilon(a)(\lambda,\theta)$. On the other hand,
if  $\lambda\in \Omega:=\{\lambda\in
P^+_{\fk}(\pi),\ \langle 
s_{\beta_l}.(\lambda)+\rho,\beta_l\rangle\in 2\N+1\}$ and
$\theta':= {(-1)}^{s_{\beta_l}.\lambda-\lambda}\theta$ 
we claim that $a$ acts on 
$\wU_{|_{1}}v_{\lambda}$ by the scalar 
$\Upsilon(a)(s_{\beta_l}.\lambda,\theta')$. 
Indeed, assume that $\lambda\in \Omega$ and 
$\theta'= {(-1)}^{s_{\beta_l}.\lambda-\lambda}\theta$.  
Then (see~\ref{qrepgen}) 
$M(q^{\lambda},\theta)\subset M(q^{s_{\beta_l}.\lambda},\theta')$. 
Moreover, if $u_{\lambda}$, $u_{s_{\beta_l}.\lambda}$ are the 
respective highest
weight vectors of these Verma modules, one has $\wU_{|_1}u_{\lambda}
\subset \wU_{|_0}u_{s_{\beta_l}.\lambda}$ and the claim follows.
Hence
$\Upsilon(a)=0$ implies $a\in\Ann_{\wU}V(q^{\lambda},\theta)$ 
for all $(\lambda,\theta)\in \Omega\times \widehat{\Gamma}$.
The density of $\Omega$ allows us to use~\ref{interannusimq} and
then to conclude.
\end{pf}

\subsection{}\label{qpfcentnotU0ev} Set 
$$\cU^o_{ev}:=\sum_{\mu\in P_{\fk}(\pi)}kK_{2\mu}\subset \cU^o.$$
The Weyl group $W$ acts on $\cU^o$ and on $\cU^o_{ev}$ in the following 
way: $$w.K_{\mu}:=q^{(\mu,w^{-1}\rho-\rho)}K_{w\mu}.$$

\begin{lem}{} \label{lemhccentwu2}
 $\Upsilon\bigl(\cC(\wU_{|_0})\bigr)\subset
  (\cU^o_{ev})^W\oplus
\xi (\cU^o_{ev})^W$
\end{lem}
\begin{pf} Firstly, we shall check that 
\begin{equation}\label{etp1q}
\Upsilon\bigl(\cC(\wU_{|_0})\bigr)\subset (\cU^o)^W\oplus \xi (\cU^o)^W.
\end{equation}
We start by fixing some notations. Recall~\ref{qnottorus}.
 For any $1\leq i\leq l$, we define
$\Gamma_i:=\{ \xi_{\mu},\ \mu\in \bigoplus_{j\not =i} (\Z\slash 2\Z)
 w_j\}$,  and the subalgebra 
$\wU_i^o:=(k\Gamma_i)\wU^o$. By definition, 
$\wU^o=\wU_i^o\oplus \xi_{w_i}\wU_i^o$.

Fix $a\in \cC(\wU_{|_0})$. For each $i=1,\ldots,l$, write  
\begin{equation}\label{eqecrq}
\Upsilon(a)=P_i+\xi_{w_i} Q_i\end{equation}
with $P_i,Q_i\in \wU^o_i$. We fix $i<l$ and show that $Q_i=0$.

Let $(\lambda,\theta)$ be in $P_{\fk}(\pi)\times\widehat{\Gamma}$ such
that  $\langle\lambda,\alpha_i\rangle
\in \N$. Consider $\theta'$  defined 
by $\theta':={(-1)}^{\langle \lambda+\rho,\alpha_i\rangle\alpha_i}\theta$.
In other words,  
$\theta'(w_j)=\theta(w_j)\ \forall j\not=i\mbox{ and } 
\theta'(w_i)={(-1)}^{\langle \lambda+\rho,\alpha_i\rangle}\theta(w_i)$.
According to~\ref{qrepgen},
$M(q^{s_{\alpha_i}.\lambda},\theta')$ is a submodule of 
$M(q^{\lambda},\theta)$. 
Moreover, as $i<l$, $\wU_{|_0}v_{s_{\alpha_i}.\lambda}\subset 
  \wU_{|_0}v_{\lambda}$ where $v_{\lambda}, v_{s_{\alpha_i}.\lambda}$ 
stand for 
the vectors of highest weight of
$M(q^{\lambda},\theta)$ and $M(q^{s_{\alpha_i}.\lambda},\theta')$.
It follows that 
\begin{equation}\label{forupsilemq}
\Upsilon(a)(s_{\alpha_i}.\lambda,{(-1)}^{\langle \lambda+\rho,\alpha_i
\rangle\alpha_i}\theta)=\Upsilon(a)(\lambda,\theta)
\end{equation} 
if $\langle \lambda,\alpha_i \rangle \in \N$.
We shall now check that formula (\ref{forupsilemq}) extends to all
$\lambda\in P_{\fk}(\pi)$. Indeed if $ \langle \lambda,\alpha_i \rangle =-1$ 
then $ \langle \lambda+\rho,\alpha_i \rangle =0$,
$s_{\alpha_i}.\lambda=\lambda$ and (\ref{forupsilemq}) is obvious.
If  $\langle \lambda,\alpha_i\rangle=-p-2,\ p\geq 0$, then 
 $s_{\alpha_i}.\lambda=\lambda+(p+1)\alpha_i$
and  $\langle s_{\alpha_i}.\lambda,
\alpha_i\rangle=p$. We can then apply (\ref{forupsilemq})  to
$s_{\alpha_i}\lambda$, which establishes  (\ref{forupsilemq}) for
$\lambda$. 

Using notation (\ref{eqecrq}),
(\ref{forupsilemq}) can be rewritten as follows
$$ P_i(s_{\alpha_i}.\lambda,\theta_i) +\theta(w_i)
{(-1)}^{\langle \lambda+\rho,\alpha_i\rangle} Q_i(s_{\alpha_i}.\lambda,
\theta_i)=P_i(\lambda,\theta_i)+\theta(w_i)Q_i(\lambda,\theta_i)$$
where $\theta_i$ is the restriction of $\theta$ to $\Gamma_i$. 
Taking successively $\theta(w_i)=\pm$ in the last equation, we obtain 
\begin{eqnarray}P_i(s_{\alpha_i}.\lambda,\theta_i)&=&P_i(\lambda,\theta_i)\\
\label{eqecrq2}
Q_i(s_{\alpha_i}.\lambda,\theta_i)&=&{(-1)}^{\langle
  \lambda+\rho,\alpha_i\rangle}
Q_i(\lambda,\theta_i)
\end{eqnarray}
Write $Q_i=\displaystyle\sum_{(\mu,\gamma)\in P_{\fg}(\pi)\times \Gamma_i}
a_{\mu,\gamma}K_{\mu}\xi_{\gamma}$. Then (\ref{eqecrq2})
implies that for all $(\lambda,\theta)\in P_{\fk}(\pi)\times 
\widehat{\Gamma_i}$,
$$ \sum_{(\mu,\gamma)\in P_{\fg}(\pi)\times \Gamma_i}
a_{\mu,\gamma}{(-1)}^{\langle \lambda,\alpha_i\rangle}
q^{(\mu,s_{\alpha_i}.\lambda)}\theta(\gamma)
+\sum_{(\mu,\gamma)\in P_{\fg}(\pi)\times \Gamma_i}
a_{\mu,\gamma}
q^{(\mu,\lambda)}\theta(\gamma)=0.$$
Since the characters $P_{\fk}(\pi)\times \widehat{\Gamma_i}\longrightarrow
k$,  $(\lambda,\theta)\mapsto 
{(-1)}^{\langle \lambda,\alpha_i\rangle}
q^{(\mu,s_{\alpha_i}.\lambda)}\theta(\gamma)$ and
$(\lambda,\theta)\mapsto 
\ q^{(\mu,\lambda)}\theta(\gamma)$ are pairwise distinct, the lemma 
of  linear independence of the characters of Dedekind forces
$Q_i=0$.

Finally, $i$ running from $1$ to $l$, we have proved that
$$\Upsilon(a)=P+\xi Q$$ where $P,Q\in \cU^o$ are invariant under the  
action of the subgroup of $W$ generated by the $s_i, i<l$.

If $\lambda$ is such that $\langle \lambda+\rho,\beta_l\rangle\in 2\N+2$,
then  $M(q^{s_{\beta_l}.\lambda},\theta)\subset M(q^{\lambda},\theta)$ with
$\wU_{|_0}v_{s_{\beta_l}.\lambda}\subset \wU_{|_0}v_{\lambda}$,
and one shows, proceeding as above, that 
\begin{equation}\label{idPQ}
P(s_{\beta_l}.\lambda)=P(\lambda),\ \ Q(s_{\beta_l}.\lambda)=Q(\lambda)
\end{equation}
for all $\lambda\in P_{\fk}(\pi)$ such that
 $\langle \lambda+\rho,\beta_l\rangle\in 2\Z$. We shall check that 
(\ref{idPQ}) actually holds for all $\lambda\in P_{\fk}(\pi)$. 
Let us treat the case of $P$ for instance. The identity (\ref{idPQ})
can be rewritten as $(s_{\beta_l}.P-P)(\displaystyle{w_l\over
 2}+\lambda')=0$ 
for all $\lambda'\in P_{\fg}(\pi)$. Write $P=\sum_{\mu\in P_{\fg}(\pi)}
a_{\mu}K_{\mu}$. Then 
$\sum_{\mu\in P_{\fg}(\pi)} q^{{1\over
    2}(\mu,w_l)}(a_{s_{\beta_l}\mu}
q^{(\mu,\rho-s_{\beta_l}\rho)}-a_{\mu})q^{(\mu,\lambda)}=0$ for
all $\lambda\in P_{\fg}(\pi)$. The linear independence of the characters 
$P_{\fg}(\pi)\rightarrow k$, $
\lambda \mapsto q^{(\mu,\lambda)}$ forces the equalities 
$a_{s_{\beta_l}\mu}
q^{(\mu,\rho-s_{\beta_l}\rho)}=a_{\mu}$ and therefore
$s_{\beta_l}.P=P$.

Finally, $P,Q$ are  $W.$-invariant and we have proved
 (\ref{etp1q}).

It remains to show that $P,Q$ are actually elements of $\cU^o_{ev}$.
For this, we should reproduce the reasoning above, analyzing now
the action of $a$ on the $M(q^{\lambda}\phi,\theta)$ where 
$\phi\in \Gamma$.

Another way to do it is to imitate~\cite{jan} 6.6, that is to consider for 
each $\phi\in \widehat{\Gamma}$, the automorphism $\sigma_{\phi}$
which keeps $\cC(\wU_{|_0})$ invariant,  and sends
 $K_{\alpha_i},E_i,F_i,\xi_i$
respectively to
$\phi(\alpha_i)K_{\alpha_i},\phi(\alpha_i)E_i,F_i,\xi_i$.
\end{pf}

\subsection{} \begin{prop}{} $\cC(\wU_{|_0})=\cA(\wU)$.
\end{prop}
\begin{pf}
By~\cite{jq} 7.17, $\Upsilon(\cZ(\cU))=(\cU^o_{ev})^W$. Since 
$\cA(\wU)\subset\cC(\wU_{|_0})$, we deduce from
lemma~\ref{lemhccentwu2} that
$\Upsilon\bigl(\cA(\wU)\bigr)=
\Upsilon(\cC(\wU_{|_0}))$. And lemma~\ref{lemhccentwu1} ends
the proof.
\end{pf}


\end{document}